A GENERAL THEORY OF (s,t)-CORE PARTITIONS

By
JOSEPH ANDREW VANDEHEY

A THESIS

Presented to the Department of Mathematics
And the Honors College at the University of Oregon
In partial fulfillment of the requirements
For the degree of
Bachelor of Science

May 2008

An Abstract of the Thesis of

| | | |
|---|---|---|
| Joseph Andrew Vandehey | for the degree of | Bachelor of Science |
| In the Department of Mathematics | to be taken | May 27, 2008 |

Title: A GENERAL THEORY OF (S,T)-CORE PARTITIONS

Approved: _________________________________________
Professor Marie Vitulli

We introduce the idea of $(s,t)$-closure and delta-sets and show that $(s,t)$-closed beta-sets which are contained set-wise in $(s,t)$-closed delta-sets are also contained partition-wise. This implies the maximal $(s,t)$-core partition theorem of Olsson and Stanton. Along the way, we reprove the $(s,s+1)$-core case, originally proved by Olsson and Stanton, using new techniques, and several lemmas and propositions regarding the containment of partitions are presented as well. Finally, we study how these apply to partitions which are $(t_1, t_2, \cdots, t_n)$-core.

# TABLE OF CONTENTS

## 1. Introduction

One of the most primitive objects of study in mathematics is the group. A group is simply a set (such as the positive rational numbers) with an operation (such as multiplication) that satisfies three rules.

First, the order of operations does not matter, so for any *a,b,c* in our group, $a(bc) = (ab)c$.

Second, there is an identity element in the group (usually called *e* or 1) such that $ea = a$ for all *a* in the group.

Third, each element has an inverse with respect to the identity. In other words, for all *a* in the group there is also a *b* in the group such that $ab = e$.

Although groups are very primitive, they can range from the extremely simplistic to the extremely complex. We can see this with a very common and popular group: the Rubik's cube. The different elements of the set are just the various different possible states for the cube, the identity element is the solved cube, and the operation is just concatenation of actions: that is, if we have two different possible states for the cube, *a* and *b*, then we can reach *a* and *b* from the solved cube by a specific set of actions, and so *ab* just means to first do all the actions needed to get from the solved cube to *a* then do all the actions needed to get from the solved cube to *b* <u>without</u> returning to the solved cube.

But while the Rubik's cube traditionally has 3x3x3 pieces, it can be extended into having 4x4x4 or 5x5x5, each time increasing the complexity of the associated group, despite the fact that all three are groups. That is because a group has little actual structure on itself just by virtue of being a group, and in the very complex cases it might be near impossible to explicitly understand what the group's structure is.

Yet we do want to understand as much as we can about a given group's structure, because groups can be found everywhere in both theoretical mathematics and applied mathematics. Group theory, then, offers us the ability to derive useful information about a group and its structure from other information we know about the group: for example, if the group has a finite number of elements and also has an element $a$ such that $a^3 = e$ then the number of elements in the group must be a multiple of 3.

Of all the various tools and techniques used by group theorists, the one most important to us in this paper is the representation. A representation of a group is a set of $n \times n$ matrices ($n \times n$ arrays of numbers with an associated multiplication rule) with a similar structure to the group itself. To be more explicit, a representation is a set of $n \times n$ matrices such that for every element $a$ in the group, there is a corresponding matrix $A$ in the representation, and if in the group $ab = c$, then in the representation $AB = C$.

To give an example, let us consider the group {-1, 1} with multiplication. One representation might look like the following.

$$1 \qquad \leftrightarrow \qquad \begin{pmatrix} 1 & 0 \\ 0 & 1 \end{pmatrix}$$

$$-1 \qquad \leftrightarrow \qquad \begin{pmatrix} -1 & 0 \\ 0 & 1 \end{pmatrix}$$

Then we can check that multiplication works the same way in the representation as it does in the original group.

$$1 \times 1 = 1 \qquad \leftrightarrow \qquad \begin{pmatrix} 1 & 0 \\ 0 & 1 \end{pmatrix} \times \begin{pmatrix} 1 & 0 \\ 0 & 1 \end{pmatrix} = \begin{pmatrix} 1 & 0 \\ 0 & 1 \end{pmatrix}$$

$$-1 \times -1 = 1 \qquad \leftrightarrow \qquad \begin{pmatrix} -1 & 0 \\ 0 & 1 \end{pmatrix} \times \begin{pmatrix} -1 & 0 \\ 0 & 1 \end{pmatrix} = \begin{pmatrix} 1 & 0 \\ 0 & 1 \end{pmatrix}$$

And so on.

But the representation need not be unique to the group. In our example, we could increase the value of $n$ and find another representation.

$$1 \leftrightarrow \begin{pmatrix} 1 & 0 & 0 \\ 0 & 1 & 0 \\ 0 & 0 & 1 \end{pmatrix} \qquad -1 \leftrightarrow \begin{pmatrix} 1 & 0 & 0 \\ 0 & -1 & 0 \\ 0 & 0 & 1 \end{pmatrix}$$

Alternately, we could keep $n$ the same and change the set of entries. For example, instead of letting the entries of the matrix be over the set of integers, we could let them be over $\mathbb{Z}_2$, where $1 + 1 = 0$, or any other field (a group-like object having two operations instead of one) of characteristic $p$ (where 1 added to itself $p$ times equals 0). While it may seem odd to write $1 + 1 = 0$, realize that in the group, 1 and 0 do not necessarily have values as we would think of them in the integers. Instead they are just the names we give to elements with very specific properties, namely 1 is the element such that $1a = a$ for all $a$, and similarly 0 is the element such that $0 + a = a$ for all $a$. We note however, that in order for our field to obey all the rules required of being a field, $p$ must be a prime.

So in $\mathbb{Z}_2$ a representation for our group might look like this.

$$1 \qquad \leftrightarrow \qquad \begin{pmatrix} 1 & 0 \\ 0 & 1 \end{pmatrix}$$

$$-1 \qquad \leftrightarrow \qquad \begin{pmatrix} 1 & 1 \\ 0 & 1 \end{pmatrix}$$

Unfortunately the representation might share the same problem that the group does in that its structure is too complex to be understood in full. The representation might also be incredibly difficult to find, but we always know that at least one representation exists (the trivial representation where every element in the group corresponds to $\begin{pmatrix} 1 & 0 \\ 0 & 1 \end{pmatrix}$). Therefore, we can use linear algebra where the properties of matrices are very well known to study the group.

In particular, we can study the trace of the matrices. The trace is just the sum of the numbers on the diagonal, so going back to our original example with

$$1 \qquad \leftrightarrow \qquad \begin{pmatrix} 1 & 0 \\ 0 & 1 \end{pmatrix}$$

$$-1 \qquad \leftrightarrow \qquad \begin{pmatrix} -1 & 0 \\ 0 & 1 \end{pmatrix}$$

The trace of the matrix on the top is $1+1=2$, and the trace of the matrix on the bottom is $-1+1=0$.

A character of a given representation then is a function $\chi$ that maps the elements of the group onto the trace of the corresponding matrix in the representation. So, in our example, $\chi(1)=2$ and $\chi(-1)=0$.

It may seem that the character is yet another layer of abstraction (and it is), but it also is one of our first breaks in the study of groups, because we do not need to know the explicit representation to calculate facts about the character of that representation. For example, if $p$ is a prime and $g$ is an element of the group such that $g^p = e$ then $\chi(g)-$ $\chi(1)$ is a multiple of $p$. We know this is true no matter what representation we base the character on.

Characters, much like integers, can be factored down into irreducible elements, and in general it is much easier to study the irreducible characters than ordinary characters. The irreducible characters over a field of characteristic $p$ are divided into so-called $p$-blocks based on how closely related the characters are (For more information, see Olsson's book [4]). We can explicitly calculate which characters are in which $p$-blocks, because to each character there is a corresponding partition, which is just an ordered list $P = (P_1, P_2, ..., P_n)$ where each $P_i$ is a positive integer and $P_i \geq P_{i+1}$. To each partition there is a way to reduce the elements of it to produce what is called a $p$-core partition. It has been proved that two characters are in the same $p$-block if and only if their corresponding partitions have the same $p$-core.

For the most part, mathematicians who study representation theory have kept $p$ fixed and studied $p$-core partitions in order to understand the distribution of characters into $p$-blocks, but a few have asked questions as to how $p$-blocks relate to $q$-blocks. Is it possible that a $p$-block equals a $q$-block? Is it possible for a $p$-block to be contained in a $q$-block? Is it possible for them to share characters at all? If so, then it implies some similarities of structure in the way the group is represented over a field of characteristic $p$ and how it is represented over a field of characteristic $q$.

In order for two characters to be present in the same $p$-block and the same $q$-block their corresponding partitions must have the same $p$-core and the same $q$-core. This has lead some mathematicians who are interested in this theory to study $(p,q)$-core partitions and to generalize them into $(s,t)$-core partitions, where $s,t$ are no longer necessarily prime.

Anderson [1] showed that there were only a finite number of $(s,t)$-core partitions and showed how to calculate them explicitly. Olsson and Stanton [5] expanded on this work giving some theorems on the relationships of $p$-blocks and $q$-blocks. They also conjectured that there exists a maximal $(s,t)$-core partition; that is, that if $P = (P_1, P_2, ..., P_n)$ is an $(s,t)$-core partition and $Q = (Q_1, Q_2, ..., Q_m)$ is the maximal $(s,t)$-core partition, then $P_i \leq Q_i$ for $i \leq n$. In their paper [5], they proved their conjecture for the $(s,s+1)$-core case.

In this paper, we will be proving Olsson and Stanton's conjecture in the general case.

In section 2, we will be providing detailed descriptions and definitions of terms we will be using in this paper as well as introducing some basic results.

In section 3, we will introduce the bead diagram as Anderson, Olsson, and Stanton have in order to provide a visual reference for the work we do. We will also introduce the idea of a delta-set and explain its significance.

In section 4, we will provide some lemmas regarding partition-containment and will use them to reprove Olsson and Stanton's result using these new lemmas.

In section 5, we will prove the general result by way of two generalizations.

In section 6, we will conclude the paper by considering partitions that are $(t_1, t_2, ..., t_n)$-core and how the various lemmas and theorems from this paper still apply in this extended case.

## 2. Preliminaries and Definitions

Of key importance to the study of (*s,t*)-core partitions are two mathematical objects: the partition and the beta-set. A partition is a finite tuple of positive integers ordered from greatest to least, while a beta-set is a finite subset of the positive integers ordered from greatest to least. The most noticeable difference between these two is that a partition is allowed to have multiple copies of the same integer, while a beta-set is not. To help differentiate the two, partitions will always be written as lists enclosed by parentheses, while beta-sets will always be written as lists enclosed by braces.

The connection between partitions and beta-sets comes from the study of hook numbers. First, the partition $P = (P_1,...,P_n)$ is represented pictorially as a series of boxes in left-aligned rows with the *i*th row having $P_i$ boxes. (We will in general write $P_i$ to refer to the *i*th element of *P*.) Each box then has a hook number. To calculate this hook number, look at the column the box is in and count the number of boxes beneath it. Then look at the row it is in and add to that the number of boxes to the right of it. Add one for the box itself and that is the hook number.

Suppose our partition is (7,6,4,4,1). Then the diagram looks like the one below.

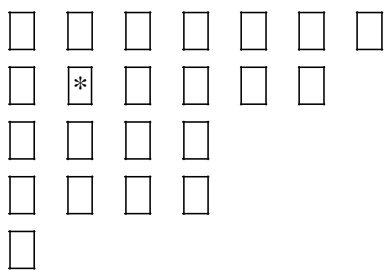

Let us calculate the hook number of the starred box by looking at its row and column and counting the required boxes.

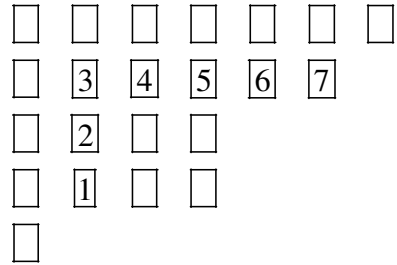

Thus the hook number for this box is 7. If we fill out the entire hook diagram for this partition, it looks like this.

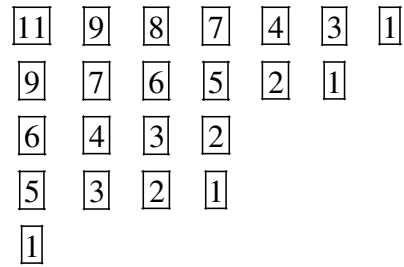

The beta-set is then just the left-most column of the hook diagram: in this case $\{11,9,6,5,1\}$.

This relationship is, in fact, a bijection between the set of all partitions and the set of all beta-sets. Given a beta-set $\beta = \{\beta_1, \beta_2, ..., \beta_n\}$, the corresponding partition is $P(\beta) = (\beta_1 - (n-1), \beta_2 - (n-2), ..., \beta_n - 0)$. Similarly, given a partition $P = (P_1, ..., P_n)$, the corresponding beta-set is $\beta(P) = \{P_1 + (n-1), P_2 + (n-2), ..., P_{n-1} + 1, P_n + 0\}$.

A partition is then said to be $t$-core if the set of all hook numbers for that partition contain no multiple of $t$. A beta-set $\beta$ is said to be $t$-core if $P(\beta)$ is $t$-core.

Similarly, a beta-set $\beta$ is said to be $t$-closed if for all $x \in \beta$ such that $x > t$, we have $x - t \in \beta$, and a partition $P$ is said to be $t$-closed if $\beta(P)$ is $t$-closed.

If *s* and *t* are relatively prime, then a partition (beta-set) is said to be (*s,t*)-core if it is both *s*-core and *t*-core and is said to be (*s,t*)-closed if it is both *s*-closed and *t*-closed. By convention we write the smaller number first, so $s < t$. Note also that whenever we refer to something being (*s,t*)-core or (*s,t*)-closed we mean that *s,t* are relatively prime.

However, more often than not we will have a beta-set that is not (*s,t*)-closed and we will want to make it (*s,t*)-closed. So in order to do that we introduce the (*s,t*)-closure of $\beta$, denote it by $\overline{\beta}$, and define it by

$$\overline{\beta} := \bigcap_{\substack{\beta' \supset \beta \\ \beta' \text{is } (s,t)-\text{closed}}} \beta' \, .$$

This definition is equivalent to $\overline{\beta} = \{x - as - bt \mid x \in \beta, a, b \geq 0, x > as + bt\}$. In the course of any given proof or formula, we fix *s* and *t*, so using the closure symbol on two different sets will mean to take the (*s,t*)-closure of both.

As a warning to the reader, other books and papers on representation theory use slightly different definitions for some of the above terms. Beta-sets, in general, do not need to be subsets just of the positive integers, but can contain 0 as well, but each such beta-set corresponds to the same partition as a beta-set that contains just positive integers, so we ignore beta-sets that contain non-positive integers as they only add an additional layer of complexity unneeded for the results of this paper. Also, it is possible to refer to the *t*-core of a partition (instead of a *t*-core partition), which is defined as the partition formed removing the *t*-hooks from the Ferrers-Young diagram of our original partition (the order of removal does not matter).

With all this in mind, we can state the first important theorem about partitions and beta-sets, which can be found in Anderson [1].

**Theorem 2.1** A beta-set $\beta$ is $(s,t)$-core if and only if $\beta$ is $(s,t)$-closed and $\beta \subset$

$\beta_{s,t} := \overline{\{st - s - t\}}$.

To prove this first we will need an extra definition and lemma. First, given a beta-set $\beta = \{\beta_1, \beta_2, ..., \beta_k\}$, let $H(\beta)$ be the set of all hook numbers (in their multiplicities) in the hook diagram of $\beta$. Then let $H_i(\beta)$ denote the subset of $H(\beta)$ containing all hook numbers in the row corresponding to $\beta_i$.

For example, if $\beta = \{6, 4, 1\}$, then the hook diagram looks like the one below.

| 6 | 4 | 3 | 1 |
| 4 | 2 | 1 |
| 1 |

In this case, $H_1(\beta) = \{6, 4, 3, 1\}$, $H_2(\beta) = \{4, 2, 1\}$, $H_3(\beta) = \{1\}$.

The following lemma will show that we can explicitly calculate $H_i(\beta)$ and hence $H(\beta)$.

**Lemma 2.2** $H_i(\beta) = \{1, 2, 3, ..., \beta_i\} \setminus \{\beta_i - \beta_j \mid j > i\}$.

**Proof of Lemma 2.2**

We work by a reverse induction, starting with $i = k$.

If we look at the actual bead diagram for $H_k(\beta)$ it looks like the diagram below.

| $\beta_k$ | | $\beta_k - 1$ | $\cdots$ | 1 |

Hence $H_k(\beta) = \{1, 2, ..., \beta_k\}$, which agrees with the statement of the lemma.

Now suppose the statement is true for a given $i$ then we wish to prove it true for $i - 1$. Let us look at the bead diagram for these two rows.

$$\boxed{\beta_{i-1}} \quad \boxed{b_1} \quad \cdots \quad \boxed{b_n} \quad \boxed{m} \quad \boxed{m-1} \quad \quad \boxed{1}$$
$$\boxed{\beta_i} \quad \boxed{a_1} \quad \cdots \quad \boxed{a_n}$$
$$\vdots \quad \vdots \quad \ddots$$

Here the various $a_i, b_i, m$ are just placeholder values we will calculate explicitly.

First, we calculate the value of $m$. Note that for $\beta_i$ to be in the box it is, we need to have that

$$\beta_i = (\text{number of boxes beneath } \boxed{\beta_i}) + (\text{number of boxes to the right of } \boxed{\beta_i}) + 1,$$

where the $+1$ at the end comes from counting the box $\boxed{\beta_i}$ itself. But then we also have that

$$\beta_{i-1} = (\text{number of boxes beneath } \boxed{\beta_{i-1}}) + (\text{number of boxes to the right of } \boxed{\beta_{i-1}}) + 1$$

$$= (\text{number of boxes beneath } \boxed{\beta_i}) + 1 + (\text{number of boxes to the right of } \boxed{\beta_{i-1}}) + 1$$

$$= (\text{number of boxes beneath } \boxed{\beta_i}) + 1 + (\text{number of boxes to the right of } \boxed{\beta_i}) + m + 1$$

$$= \beta_i + m + 1.$$

By the same logic, we also have that $b_j = a_j + m + 1 = a_j + \beta_{i-1} - \beta_i$.

Thus we can generate $H_{i-1}(\beta)$ from $H_i(\beta)$.

$$H_{i-1}(\beta) = \left( H_i(\beta) + \beta_{i-1} - \beta_i \right) \cup \{1, 2, ..., m\}$$

$$= \left( \left( \{1, 2, 3, ..., \beta_i\} \setminus \{\beta_i - \beta_j \mid j > i\} \right) + \beta_{i-1} - \beta_i \right) \cup \{1, 2, ..., \beta_{i-1} - \beta_i - 1\}$$

$$= \left( \left( \{1, 2, 3, ..., \beta_i\} + \beta_{i-1} - \beta_i \right) \setminus \left( \{\beta_i - \beta_j \mid j > i\} + \beta_{i-1} - \beta_i \right) \right) \cup \{1, 2, ..., \beta_{i-1} - \beta_i - 1\}$$

$$= \left( \{1 + \beta_{i-1} - \beta_i, 2 + \beta_{i-1} - \beta_i, ..., \beta_{i-1}\} \setminus \{\beta_{i-1} - \beta_j \mid j > i\} \right) \cup \{1, 2, ..., \beta_{i-1} - \beta_i - 1\}$$

$$= \{1, 2, ..., \beta_{i-1}\} \setminus \{\beta_{i-1} - \beta_j \mid j > i - 1\}$$

As desired. Thus, by the induction hypothesis, $H_i(\beta) = \{1, 2, 3, ..., \beta_i\} \setminus \{\beta_i - \beta_j \mid j > i\}$ in general.

Q.E.D.

Now to prove the theorem, remember that $\beta$ being $t$-core is equivalent to $tk \notin H(\beta)$ for any positive integer $k$.

We can rewrite this as $\beta$ is not $t$-core if and only if there exist $i, k$ such that $tk \in H_i(\beta)$. But $tk \in H_i(\beta)$ if and only if there does not exist a $j$ such that $tk = \beta_i - \beta_j$ or $\beta_i = tk$.

Putting this together we see that $\beta$ is $t$-core if and only if given $\beta_i - tk > 0$ we have some $j$ for which $\beta_i - tk = \beta_j$, and $\beta_i \neq tk$. But saying that $\beta_i - tk > 0$ implies that we have some $j$ for which $\beta_i - tk = \beta_j$ is equivalent to saying that the entire sequence $\beta_i$, $\beta_i - t, \beta_i - 2t, ..., \beta_i - kt$ must be in $\beta$, i.e. that $\beta$ is $t$-closed.

So in order for $\beta$ to be $(s,t)$-core, it must be $(s,t)$-closed and contain no members of the form $as + bt$, where $a, b$ are positive. This statement is equivalent to Theorem 2.1, since the numbers which cannot be represented as $as + bt$ are precisely the numbers in $\beta_{s,t}$. So we have proved Theorem 2.1.

This theorem implies a number of useful statements, including the fact that there are only finitely many $(s,t)$-core partitions. It also implies a simple test for whether a given partition is $(s,t)$-core, without needing to calculate the entire hook diagram: namely, just compute the beta-set of that partition and then check the two conditions above. We can rewrite this as a corollary.

**Corollary 2.3** A partition $P$ is $(s,t)$-core if and only if the following conditions are satisfied for all $x \in \beta(P)$:

$(i)$ $\quad x < s$ or $x - s \in \beta(P)$;

$(ii)$ $\quad x < t$ or $x - t \in \beta(P)$;

$(iii)$ $\quad x = st - as - bt$, for some positive integers $a, b$.

At times it is also useful to describe a beta-set using as few elements as possible. For that we use a generating set, which, for an $(s,t)$-closed beta-set $\beta$, is a subset $\alpha \subset \beta$ such that $\bar{\alpha} = \beta$. The minimal generating set for $\beta$ is the smallest such generating set and can be written as the intersection of all generating sets for $\beta$. We generally refer to the elements of the minimal generating set as the generators of $\beta$. $n$ is a generator of $\beta$ if and only if $n \in \beta$, but $n + s \notin \beta$ and $n + t \notin \beta$.

We can add a positive integer $a$ to a partition (or beta-set) element-wise and define this operation simply by $(P_1, P_2, ..., P_n) + a = (P_1 + a, P_2 + a, ..., P_n + a)$. We define subtraction by a positive integer similarly, except with the added rule that if this operation causes an element to shrink below 1, that element is deleted, because we do not allow non-positive numbers in partitions or beta-sets. Because of this additional rule, we note that $P + a - a$ does not necessarily equal $P - a + a$.

Finally, we need some way to compare two partitions or beta-sets. So, given two partitions $P = (P_1, ..., P_n), Q = (Q_1, ..., Q_m)$, we write $P < Q$ (read $P$ is contained partition-wise in $Q$) if $n \leq m$ and for all $i \leq n$, $P_i \leq Q_i$. Given two beta-sets $\beta = (\beta_1, ..., \beta_n)$, $\beta' = (\beta'_1, ..., \beta'_m)$, we write $\beta \prec \beta'$ (again read $\beta$ is contained partition-wise in $\beta'$) if $P(\beta) < P(\beta')$; that is, $\beta \prec \beta'$ if $n \leq m$ and for all $i \leq n$, $\beta_i \leq \beta'_i - (m - n)$.

The notation $\beta \subset \gamma$ can also be read as $\beta$ is contained in $\gamma$, so when there is confusion we refer to this kind as set-wise containment and the notation defined in the previous paragraph as partition-wise containment.

We can now present the theorem we will be proving over the course of this paper:

**Theorem 2.4 (Maximal Theorem)**

There exists a maximal $(s,t)$-core partition under partition-wise containment. In particular, if $\beta$ is any $(s,t)$-core beta-set, then $\beta \prec \beta_{s,t} = \overline{\{st - s - t\}}$.

But in order to understand this theorem better and to begin proving it, we need to understand some basic rules of partition inclusions.

**Proposition 2.5** Let $\beta$ and $\gamma$ be beta-sets. If every element in $\beta$ is larger than every element in $\gamma$, then $P(\beta \cup \gamma) = P(\beta - |\gamma|) \cup P(\gamma)$, where the union symbol on the right hand side denotes concatenation of partitions.

**Proof of Proposition 2.5**

We simply do a direct comparison of both sides. Given $\beta = \{\beta_1, \beta_2, ..., \beta_n\}$, $\gamma = \{\gamma_1, \gamma_2, ..., \gamma_m\}$, we have that

$$P(\beta - |\gamma|) \cup P(\gamma)$$

$$= (\beta_1 - m - (n-1), \beta_2 - m - (n-2), ..., \beta_n - m) \cup (\gamma_1 - (m-1), \gamma_2 - (m-2), ..., \gamma_m)$$

$$= (\beta_1 - (m+n-1), \beta_2 - (m+n-2), ..., \beta_n - m, \gamma_1 - (m-1), \gamma_2 - (m-2), ..., \gamma_m)$$

$$= P(\beta \cup \gamma).$$

$$\text{Q.E.D.}$$

This first proposition allows us to study the properties of a given partition by breaking it into two smaller partitions; this will be particularly useful in Proposition 2.8. First, though, we show a simple case of partition containment.

**Proposition 2.6** If $\beta = \{\beta_1, \beta_2, ..., \beta_n\}$, $\gamma = \{\gamma_1, \gamma_2, ..., \gamma_m\}$ are two beta-sets such that $|\beta| = |\gamma|$ and $\beta_i \leq \gamma_i$ for all $i \leq n$, then $\beta \prec \gamma$.

**Proof of Proposition 2.6**

For all $i \leq n$, $P(\beta)_i = \beta_i - (n-i) \leq \gamma_i - (n-i) = P(\gamma)_i$, so $P(\beta) < P(\gamma)$.

Q.E.D.

**Corollary 2.7** Let $\beta = \{\beta_1, \beta_2, ..., \beta_n\}$ be a beta-set. If $\beta_j + a \notin \beta$ for some $a > 0$ then $\beta \prec (\beta \cup \{\beta_j + a\}) \setminus \{\beta_j\}$.

**Proof of Corollary 2.7**

Let $\gamma = (\beta \cup \{\beta_j + a\}) \setminus \{\beta_j\}$, then $\gamma = \{\beta_1, \beta_2, ..., \beta_{k-1}, \beta_j + a, \beta_k, \beta_{k+1}, ..., \beta_{j-1}, \beta_{j+1}, \beta_{j+2}, ..., \beta_n\}$, where $\beta_{k-1} > \beta_j + a > \beta_k$.

We then apply Proposition 2.6. Clearly $|\beta| = |\gamma|$, moreover we have that $\beta_i = \gamma_i$ for $i < k$ and $i > j$. We also have that $\gamma_k = \beta_j + a > \beta_k$. Finally for all other $i$, $\gamma_i = \beta_{i-1} > \beta_i$, since beta-sets are strictly decreasing. Thus $\beta_i \leq \gamma_i$ for all $i \leq n$ and so $\beta \prec \gamma$.

Q.E.D.

Of course, we can apply Corollary 2.7 as many times as we wish, increasing the value of as many different elements as we like.

Returning to the idea of Proposition 2.5 , we want to not only split partitions into smaller elements but then to use this as a way to show partition containment.

**Proposition 2.8** Given two beta-sets, $\beta = \{\beta_1, \beta_2, ..., \beta_n\}$, $\gamma = \{\gamma_1, \gamma_2, ..., \gamma_m\}$, we have that $\beta \prec \gamma$ if, and only if, $n \leq m$ and for every integer $i \leq n$, both the following statements are true:

(i) $\{\beta_1, \beta_2, ..., \beta_i\} - (|\beta| - i) \prec \{\gamma_1, \gamma_2, ..., \gamma_i\} - (|\gamma| - i)$; and,

(ii) $\{\beta_{i+1}, \beta_{i+2}, ..., \beta_n\} \prec \{\gamma_{i+1}, \gamma_{i+2}, ..., \gamma_m\}$.

Alternately, given two partitions $A = (A_1, A_2, ..., A_n), B = (B_1, B_2, ..., B_m)$, we have that $A \leq B$ if, and only if, $n \leq m$ and for every integer $i \leq n$, both the following statements are true:

(i) $(A_1, A_2, ..., A_i) \leq (B_1, B_2, ..., B_i)$; and,

(ii) $(A_{i+1}, A_{i+2}, ..., A_n) \leq (B_{i+1}, B_{i+2}, ..., B_m)$.

**Proof of Proposition 2.8**

We prove the second half first. Note that statement (i) says that $A_j \leq B_j$ for $1 \leq j \leq i$, statement (ii) says that $A_j \leq B_j$ for $i+1 \leq j \leq n$, and we also have that $n \leq m$. This is equivalent to saying that $A_j \leq B_j$ for $1 \leq j \leq n$ and $n \leq m$, which is itself equivalent to saying that $A \leq B$.

For the first half, note that statement (i) says that $\beta_j - (i - j) - (n - i) \leq \gamma_j - (i - j) - (m - i)$ for $1 \leq j \leq i$, and statement (ii) says that $\beta_j - (n - j) \leq \gamma_j - (m - j)$ for $i+1 \le j \le n$, and we also have that $n \le m$. This is equivalent to saying that $\beta_j - (n - j) \le \gamma_j - (m - j)$ for $1 \le j \le n$ and $n \le m$, which is itself equivalent to saying that $\beta \prec \gamma$.

<div align="right">Q.E.D.</div>

**Proposition 2.9** A beta-set $\beta$ is $(s,t)$-closed if and only if for all natural numbers $y \notin \beta$, we also have that $y + s, y + t \notin \beta$.

**Proof of Proposition 2.9**

We already know that $\beta$ is $(s,t)$-closed if and only if for all $x \in \beta$, if $x > s$, then $x - s \in \beta$, and if $x > t$, then $x - t \in \beta$.

First we show that if $\beta$ is $(s,t)$-closed then for all natural numbers $y \notin \beta$ we also have $y + s, y + t \notin \beta$. Suppose we have a natural number $y \notin \beta$, but $y + s \in \beta$. Since $y + s > s$ and $\beta$ is $(s,t)$-closed, $(y + s) - s = y$ must be an element of $\beta$, which is a contradiction. In the same way, if $y \notin \beta$, but $y + t \in \beta$, we also reach a contradiction.

Now suppose we know that for all natural numbers $y \notin \beta$, we have that $y + s$, $y + t \notin \beta$, and we want to show this would imply that $\beta$ is $(s,t)$-closed. Suppose $x \in \beta$ and $x > s$ but $x - s \notin \beta$. Then if we set $y = x - s$, we see that $y + s = (x - s) + s = x \notin \beta$ which is a contradiction. Similarly $x \in \beta$ and $x > t$ but $x - t \notin \beta$ yields a contradiction.

<div align="right">Q.E.D.</div>

These propositions are somewhat trivial in that their proofs are brute-force checks of partition containment; however, by proving them, we do not need to perform the same brute-force checks later when proving more complex lemmas and theorems.

## 3.  Visual representations and delta-sets

We can represent (*s*,*t*)-core beta-sets pictorially by using a bead diagram (for more information about bead diagrams and why their construction is useful, see Anderson [1]).  To create the bead diagram we take the elements of $\beta_{s,t}$ and arrange them like so:

$$
\begin{array}{cccc}
\vdots & & & \\
st-3s-t & \cdot\cdot^{\cdot} & & \\
st-2s-t & st-2s-2t & \cdot\cdot^{\cdot} & \\
st-s-t & st-s-2t & st-s-3t & \cdots
\end{array}
$$

So that above each number is a number *s* less than it, and to the right of each number is a number *t* less than it.  As beta-sets contain only positive numbers, the bead diagram is finite.  Since *s* is smaller than *t* by convention, the bead diagram is typically taller than it is wide.  As an example, for (5,7)-core the bead diagram looks like:

$$
\begin{array}{cccc}
3 & & & \\
8 & 1 & & \\
13 & 6 & & \\
18 & 11 & 4 & \\
23 & 16 & 9 & 2
\end{array}
$$

Then to represent a given beta-set on the bead diagram we circle the numbers in the diagram which are elements of the beta-set.  This allows us to further simplify the test for (*s*,*t*)-core given in Corollary 2.3: every circled number should have a circled number above and to the right of it (conditions (*i*) and (*ii*)) and every number in the beta-set should appear in the diagram (condition (*iii*)).

However, because we will often need to consider (*s*,*t*)-closed beta-sets which are not necessarily (*s*,*t*)-core we introduce the notion of the extended bead diagram, which we create by adding an infinite number of new rows beneath the diagram by adding *s* and an infinite number of new columns to the left by adding *t*. To continue the (5,7)-core bead diagram out a bit would give:

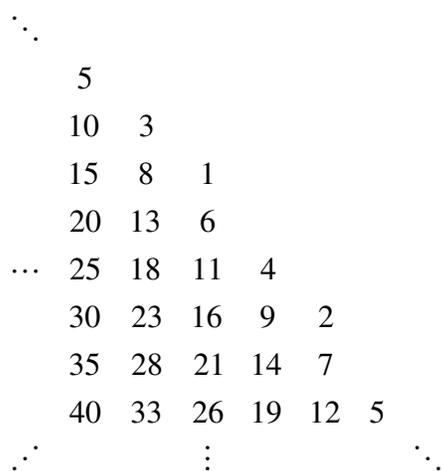

It should be noted that every positive integer occurs an infinite number of times in the extended bead diagram (since starting at a given number, going *s* columns to the left and *t* rows up gives the same number again). So we must extend also the way we think about beta-sets being represented in the bead diagram: that is, we must circle every appearance of a number from the beta-set, despite there only being one copy of that number in the beta-set itself.

Note that it is possible for delta-set to become so large that it eventually overlaps itself. For example, in the delta-set $\overline{\{40\}}$ in the extended bead diagram above, we could reach the element 5 either by starting at the generator and moving 7 rows up, or by moving 5 columns to the right. In general, we will treat these overlapping delta-sets like any other delta-sets, but we will point out in the paper when they need to be treated with separate cases.

We denote the height and width of a beta-set $\beta$ by the symbols $h(\beta)$ and $w(\beta)$ where the height is the number of residue classes left when reducing $\beta$ modulo $t$, and the width the number left when reducing modulo $s$. If we restrict our attention to the bead diagram, we notice that all the elements in a given column have the same residue class modulo $s$, and all the elements in a given row have the same residue class modulo $t$. Thus the height is the number of distinct rows containing elements from $\beta$ and the width is the number of distinct columns containing elements from $\beta$. We say distinct here, because like elements, rows and columns repeat themselves an infinite number of times in the extended bead diagram.

In general we will constrain ourselves to using the regular bead diagram and only add in extra rows and columns as needed. Using the bead diagram we can visualize several different types of beta-sets.

Beta-sets that are not $(s,t)$-closed:

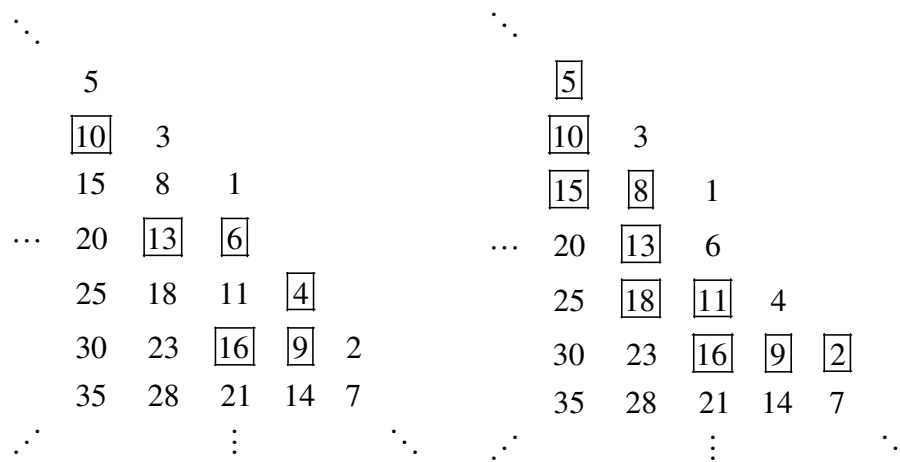

Beta-sets that are $(s,t)$-closed but not $(s,t)$-core:

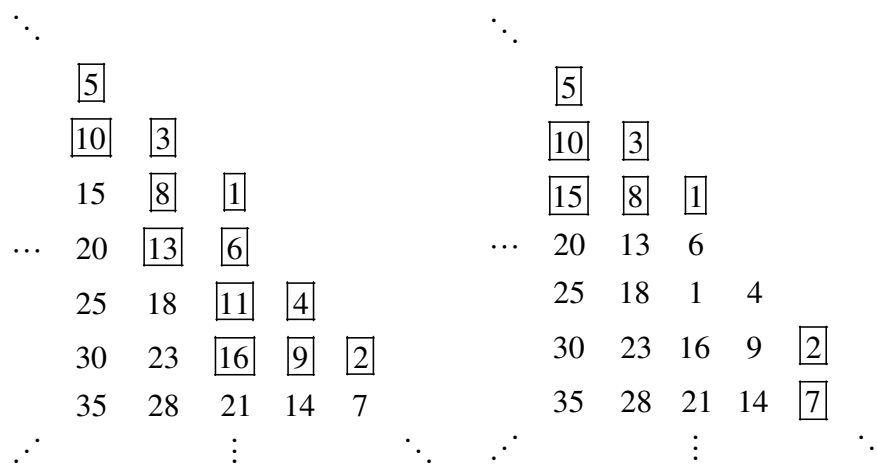

Beta-sets that are $(s,t)$-core:

```
 [3]                        3
 [8] [1]                    8  [1]
 13  [6]                   13  [6]
 18  11  4                 18 [11] [4]
 23  16  9  [2]            23  16 [9] [2]
```

And finally a special class of $(s,t)$-closed beta-sets:

```
  ⋰
      [5]
     [10] [3]                        3
     [15] [8] [1]                    8  [1]
 ··· 20   13  6                     13  [6]
     25   18  11  4                 18 [11] [4]
     30   23  16  9   2             23  16  9   2
     35   28  21  14  7
 ⋰          ⋮           ⋱
```

These beta-sets are $(s,t)$-closed and have only one generator. We call them delta-sets. These delta-sets have a number of nice properties which are summed up in the following theorem.

**Theorem 3.1** Given $i, j \in \mathbb{N}$, and $s,t$ relatively prime, the following hold:

$(i)$  $\qquad \overline{\{i\}} - j = \overline{\{i - j\}}$ ;

$(ii)$ $\qquad \overline{\overline{\{i\}}+j} = \overline{\{i+j\}}$ ;

$(iii)$ $\qquad$ If $\{s,t\} \cap \overline{\{i\}} = \varnothing$ , then $\left|\overline{\{i+1\}}\right| - \left|\overline{\{i\}}\right| = 0$ ; otherwise $\left|\overline{\{i+1\}}\right| - \left|\overline{\{i\}}\right| = 1$ ;

$(iv)$ $\qquad$ For any $n \geq 0$ , $\exists i$ s.t. $\left|\overline{\{i\}}\right| = n$ ;

$(v)$ $\qquad$ We have $a < b$ if and only $\overline{\{a\}} \prec \overline{\{b\}}$ ;

$(vi)$ $\qquad$ $\left|\overline{\{a\}}\right| < \left|\overline{\{b\}}\right|$ implies $a < b$ .

## Proof of Theorem 3.1

$(i)$ $\qquad$ If $j \geq i$ , then both $\overline{\{i\}}-j$ and $\overline{\{i-j\}}$ are empty. If $j < i$ , then

$$\overline{\{i\}} - j = \{i - as - bt \mid a,b \geq 0 \text{ and } as+bt < i\} - j$$

$$= \{i - j - as - bt \mid a,b \geq 0 \text{ and } as+bt < i-j\}$$

$$= \overline{\{i-j\}} .$$

The second equality holds because beta-sets contain no non-positive integers.

$(ii)$ $\qquad$ $\overline{\overline{\{i\}}+j} = \overline{\{i - as - bt \mid a,b \geq 0 \text{ and } as+bt < i\} + j}$

$\qquad = \overline{\{i + j - as - bt \mid a,b \geq 0 \text{ and } as+bt < i\}}$

$\qquad = \overline{\{i + j - a's - b't \mid a',b' \geq 0 \text{ and } a's+b't < i+j\}}$

$\qquad = \overline{\{i+j\}}$

$(iii)$ $\qquad$ $\overline{\{i+1\}} = \{i + 1 - as - bt \mid a,b \geq 0 \text{ and } as+bt < i+1\}$

Suppose for some integer $j$, $j \in \overline{\{i+1\}}$. This says that there exist $a, b \geq 0$ such that $i + 1 - as - bt = j$. By simple algebra, this means that there exists $a, b \geq 0$ such that $i - as - bt = j - 1$. But since $\overline{\{i\}} = \{i - as - bt \mid a, b \geq 0 \text{ and } as + bt < i\}$, either $j - 1 \in \overline{\{i\}}$ or $j - 1 \leq 0$, but the latter case can only happen if $j = 1$.

By the same logic, if $j \in \overline{\{i\}}$ then $j + 1 \in \overline{\{i+1\}}$.

From this we can tell that $0 \leq \left| \overline{\{i+1\}} \right| - \left| \overline{\{i\}} \right| \leq 1$, and if the right hand equality holds then $1 \in \overline{\{i+1\}}$.

Now suppose $1 \in \overline{\{i+1\}}$, then since the set in question is a delta-set and can have only one generator, either 1 is a generator, in which case $i = 0$ which is impossible, or 1 is not a generator and hence $1 + s \in \overline{\{i+1\}}$ or $1 + t \in \overline{\{i+1\}}$. But as was shown earlier this means that then $s \in \overline{\{i\}}$ or $t \in \overline{\{i\}}$. So $\left| \overline{\{i+1\}} \right| - \left| \overline{\{i\}} \right| = 1$ in this case.

If instead $s, t \notin \overline{\{i\}}$ then $s + 1, t + 1 \notin \overline{\{i+1\}}$, so $1 \notin \overline{\{i+1\}}$. Thus $\left| \overline{\{i+1\}} \right| - \left| \overline{\{i\}} \right| = 0$ in this case.

($iv$)     We know that the delta-set $\overline{\{1\}} = \{1\}$ has only one element. Now we just compare each delta-set with each successive delta-set (that is, add 1 to the generator by part ($ii$)). By part ($iii$), we know that the difference in size of two successive delta-sets can at most be 1, so it suffices to show that no matter how many elements a delta-set has, by adding enough value to the generator, you can get a delta-set with at least one more element.

So, suppose by contradiction that the size of delta-sets stops at $n$. Let $i$ be the smallest value for which $\left|\overline{\{i\}}\right| = n$. By our discussion in the proof of part ($iii$) we know that $1 \in \overline{\{i\}}$. Therefore $s \in \overline{\{i+s-1\}}$, but again by our discussion of the proof of part ($iii$), this means that $1 \in \overline{\{i+s\}}$, and hence $\left|\overline{\{i+s\}}\right| > \left|\overline{\{i\}}\right|$.

($v$)    To prove that $a < b$ implies $\overline{\{a\}} \prec \overline{\{b\}}$, it is enough to show that $\overline{\{a\}} \prec \overline{\{a+1\}}$.

Suppose $\overline{\{a\}} = \{a_1, a_2, ..., a_n\}$, then by the proof of part ($iii$) we have that $\overline{\{a+1\}} = \{a_1+1, a_2+1, ..., a_n+1\} \cup A$, where $A \subset \{1\}$.

Then for $i \leq n$, $P(\overline{\{a+1\}})_i = (a_i+1) - (n-i) - |A|$. But since $|A|$ is at most 1, we have that $P(\overline{\{a+1\}})_i \geq a_i - (n-i) = P(\overline{\{a\}})_i$. Thus $\overline{\{a\}} \prec \overline{\{a+1\}}$.

To prove the reverse implication it is enough to show that $\overline{\{a+1\}} \not\prec \overline{\{a\}}$. Let $P(\overline{\{a\}}) = (P_1, P_2, ..., P_n)$, then if $A = \varnothing$, $P(\overline{\{a+1\}}) = (P_1+1, P_2+1, ..., P_n+1)$, and if $A = \{1\}$, then $P(\overline{\{a+1\}}) = (P_1, P_2, ..., P_n, 1)$. By direct comparison, we have that $P(\overline{\{a+1\}}) \not\prec P(\overline{\{a\}})$.

($vi$)    Suppose $\left|\overline{\{a\}}\right| < \left|\overline{\{b\}}\right|$.

If $a = b$, then $\overline{\{a\}} = \overline{\{b\}}$, so $\left|\overline{\{a\}}\right| = \left|\overline{\{b\}}\right|$, which contradicts our original assumption.

If $a > b$, then we can consider the chain of delta-sets $\overline{\{b\}} \prec \overline{\{b+1\}} \prec \overline{\{b+2\}} \prec \cdots \prec \overline{\{a\}}$ (the partition-wise containment holds by part ($v$)). By part ($iii$), $\left|\overline{\{b\}}\right| \leq \left|\overline{\{b+1\}}\right| \leq$ $\left|\overline{\{b+2\}}\right| \le \cdots \le \left|\overline{\{a\}}\right|$, which importantly shows that $\left|\overline{\{b\}}\right| \le \left|\overline{\{a\}}\right|$, a contradiction to our original assumption.

The only remaining possibility is that $a < b$.

Q.E.D.

$\beta_{s,t}$, from the maximal theorem, is itself a delta-set and, in fact, the maximal theorem is true not because of any special property unique to $\beta_{s,t}$, but due to a special property that all delta-sets possess. This special property gives us our first generalization of the maximal theorem.

### Theorem 3.2 (First generalization of the maximal theorem)

Let $\Delta$ be any (*s,t*)-closed delta-set and $\beta$ be any (*s,t*)-closed beta-set such that $\beta \subset \Delta$. Then $\beta \prec \Delta$.

We will prove this statement in section 5. Since $\beta_{s,t}$ is itself a delta-set, the first generalization trivially proves the maximal theorem.

## 4. Reproving the (*s,s*+1) case

We begin our work towards proving the maximal theorem by studying the (*s,s*+1)-core case that Olsson and Stanton originally proved. We will use the construction of canonical forms of beta-sets as they did, but will also prove new lemmas regarding partition containment. We will then show the strength of these new lemmas by using them to simplify the end of the Olsson and Stanton proof.

**Lemma 4.1** For any (*s,t*)-closed beta-set $\beta$ and any positive integer $k$, $\beta \prec \overline{\beta + k}$.

## Proof of Lemma 4.1

By definition, $\overline{\beta + k}$ contains $\beta + k$ set-wise. So $\overline{\beta + k} = (\beta + k) \cup A$, where $A$ is some (possibly empty) subset of the positive integers disjoint from $\beta + k$; however, we can be more specific. $A$ consists of integers $x \notin \beta + k$ where $x + as \in \beta + k$ or $x + bt \in \beta + k$ for some positive $a,b$. This implies as well that $x - k \notin \beta$ and either $x + as - k \in \beta$ or $x + bt - k \in \beta$. But $\beta$ is (*s,t*)-closed, so this is only possible if $x \leq k$ (otherwise we would violate Proposition 2.9). Thus $A$ is some subset of $\{1, 2, ..., k\}$.

Now we can do a direct comparison of the partitions. If $\beta = \{\beta_1, \beta_2, ..., \beta_n\}$ then since all the elements in $A$ are smaller than all the elements in $\beta + k$, we have also that

$$P(\overline{\beta + k}) = P((\beta + k) \cup A)$$

$$= (P(\beta + k) - |A|) \cup P(A)$$

$$= ((\beta_1 + k - (n-1), \beta_2 + k - (n-2), ..., \beta_{n-1} + k - 1, \beta_n + k) - |A|) \cup P(A)$$

$$= ((\beta_1 - (n-1), \beta_2 - (n-2), ..., \beta_{n-1} - 1, \beta_n) + k - |A|) \cup P(A)$$

$$= (P(\beta) + k - |A|) \cup P(A).$$

But from what we discovered about $A$, $k - |A| \geq 0$, so an element-by-element comparison shows $P(\beta) < (P(\beta) + k - |A|) < P(\overline{\beta + k})$. Thus $\beta \prec \overline{\beta + k}$.

<div align="right">Q.E.D.</div>

We can visualize Lemma 4.1 on the extended bead diagram by starting with a beta-set then creating a second beta-set by rigidly moving the first beta-set into a new position and filling in additional elements as necessary to make it $(s,t)$-closed. Then we know that the first beta-set is contained partition-wise in the second one.

For example, if we consider $(5,7)$-core partitions, then suppose $\beta = \{6,1\}$ and

$k = 12$. Then pictorially we have

```
        ⋱
          5
         10   3
         15   8   1
β = ...  20   13  6
         25   18  11   4
         30   23  16   9   2
         35   28  21   14  7
       ⋰          ⋮           ⋱
```

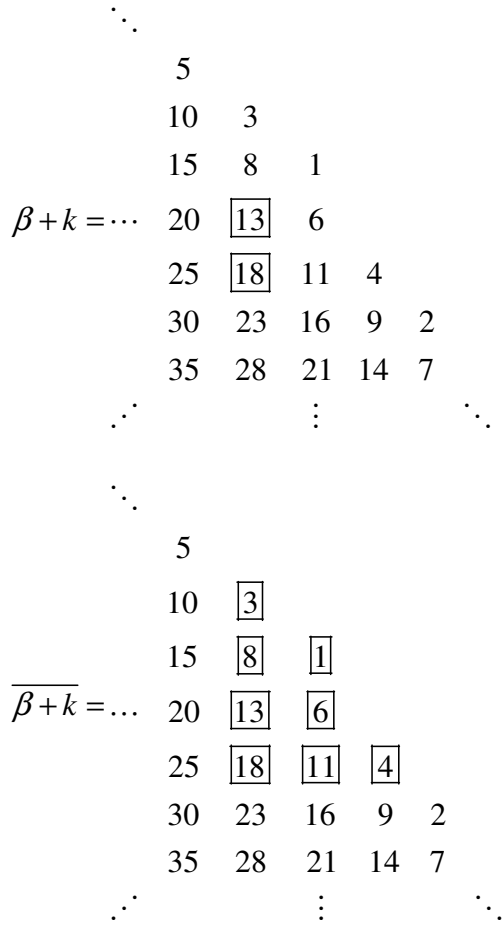

We can then check to see that $P(\beta) = (5,1) < (11,7,6,4,3,2,2,1) = P(\overline{\beta+k})$ as desired.

If we restrict ourselves to $(s,t)$-core beta-sets then we have the following useful corollary.

**Corollary to Lemma 4.1** If $\beta$ is an $(s,t)$-core beta-set that contains no elements from the bottom row of the bead diagram, $\beta \prec \overline{\beta+s}$ and $\overline{\beta+s}$ is $(s,t)$-core. Similarly if $\beta$ contains no elements from the leftmost column of the bead diagram, $\beta \prec \overline{\beta+t}$ and $\overline{\beta+t}$ is $(s,t)$-core.

Using Lemma 4.1, we can easily reprove part ($v$) of Theorem 3.1. If $a < b$ then $\overline{\{a\}} \prec \overline{\overline{\{a\}} + (b-a)} = \overline{\{b\}}$, so the delta-sets form a totally ordered set under partition-wise containment. Therefore it makes sense to talk about the largest element of a finite collection of delta-sets and we can now state our second generalization of the maximal theorem.

**Theorem 4.2 (Second generalization of the maximal theorem)**

Let $M\Delta_n$ be the largest (under partition containment) ($s,t$)-closed delta-set containing $n$ elements. Then for any ($s,t$)-closed beta-set $\beta$ containing $n$ elements, $\beta \prec \overline{M\Delta_n + s}$.

We will also prove this in section 5.

It is not true in general that $\beta \prec M\Delta_n$. As we check element by element, the containment will hold for a while but in the last few elements it breaks down. We can see roughly why this will happen: consider the sets $\widetilde{M\Delta_n} = M\Delta_n \cap \{1, 2, \cdots, s\}$ and $\tilde{\beta} = \beta \cap \{1, 2, \cdots, s\}$, which, in the bead diagram are the elements on the top of each column. We have only one possibility for the elements in $\widetilde{M\Delta_n}$, and they may not be especially large values, but we have many possibilities for $\tilde{\beta}$, some of which might contain large enough values that in a direct comparison they will overshadow the terms in $\widetilde{M\Delta_n}$.

The following lemma has much more esoteric requirements than our previous one, but we can sometimes force beta-sets to match its requirements in order to find more interesting relationships. In fact we will use the previous lemma to force beta-sets to match the requirements of this lemma.

**Lemma 4.3** Given two beta-sets $\beta, \gamma$, such that $\beta \subset \gamma$ and $\beta \prec \gamma$, we also have that $\beta \cup \{k\} \prec \gamma$ where $k = \min(\gamma \setminus \beta)$. (Note that this lemma does not require that the beta-sets even be $(s,t)$-closed.)

**Proof of Lemma 4.3**

Let $\beta = \{\beta_1, \beta_2, ..., \beta_n\}$.

Then $\beta \cup \{k\} = \{\beta_1, \beta_2, ..., \beta_i, k, \beta_{i+1}, ..., \beta_n\}$ for some $i$.

Since $\beta \subset \gamma$ and $k$ is the smallest element of their difference, $k$ must be an element of $\gamma$. Thus, we have $\gamma = \{\gamma_1, \gamma_2, ..., \gamma_j, \gamma_{j+1}, ..., \gamma_m\} = \{\gamma_1, \gamma_2, ..., \gamma_j, k, \beta_{i+1}, ..., \beta_n\}$. Here, $i$ must be less than $j$ because otherwise $\beta$ would have more elements than $\gamma$.

We know by our assumptions that $P(\beta)_x \leq P(\gamma)_x$ for all $x \leq n$, and in particular, $P(\beta \cup \{k\})_x = P(\beta)_x - 1 < P(\gamma)_x$ for $x \leq i$.

Thus by Proposition 3, in order to show that $\beta \cup \{k\} \prec \gamma$ we need only show that $P(\beta \cup \{k\})_x \leq P(\gamma)_x$ for $i < x \leq n+1$.

However, these remaining elements are precisely $\{k, \beta_{i+1}, ..., \beta_n\} = \{\gamma_{j+1}, \gamma_{j+2}, ..., \gamma_m\}$ in $\beta \cup \{k\}$. Thus, $P(\beta \cup \{k\})_x = \gamma_{j+(x-i)} - (m - j - (x-i))$ for $i < x \leq n+1$.

We already know that since beta-sets do not repeat elements and are arranged in descending order, $\gamma_x \geq \gamma_{x+1} + 1$, so repeating this several times we get $\gamma_x = \gamma_{i+(x-i)}$

$\geq \gamma_{j+(x-i)} + j - i$, which in turn gives $P(\beta \cup \{k\})_x = \gamma_{j+(x-i)} - (m - j - (x - i))$

$\leq \gamma_x - (m - x) = P(\gamma)_x$ for $n + 1 \geq x > i$ as desired. Thus $\beta \cup \{k\} \prec \gamma$.

<div align="right">Q.E.D.</div>

**Corollary to Lemma 4.3** If $\beta, \gamma$ are as in Lemma 4.3 but are now $(s,t)$-closed (respectively, $(s,t)$-core), then $\beta \cup \{k\}$ is $(s,t)$-closed (respectively, $(s,t)$-core).

**Proof of Corollary to Lemma 4.3**

This is trivially true if we consider $\beta$ and $\gamma$ to be $(s,t)$-closed, because if $k - s$ or $k - t$ were not an element of $\beta$, then $k$ would not be the smallest element of the difference between $\gamma$ and $\beta$. For $\beta$ and $\gamma$ $(s,t)$-core, note that $\beta$ and $\gamma$ are contained set-wise in $\beta_{s,t}$ and both are $(s,t)$-closed. We have that $\beta \cup \{k\}$ is $(s,t)$-closed by the first half of this proof, and since $\beta \cup \{k\} \subset \gamma$, $\beta \cup \{k\}$ satisfies both conditions for being $(s,t)$-core as well.

<div align="right">Q.E.D.</div>

Now we will give an example of a practical use of Lemmas 4.1 and 4.3 which will be used in the proof of the $(s, s+1)$-core case. Assume $\beta$ is an $(s,t)$-core beta-set with no elements in the bottom row of the bead diagram, then by the Corollary to Lemma 4.1, $\overline{\beta + s}$ is $(s,t)$-core as well, but more importantly, we have both $\beta \subset \overline{\beta + s}$ and

$\beta \prec \overline{\beta + s}$. Now we can apply Lemma 4.3 multiple times using the set $\{k_1, k_2, ..., k_n\} =$

$\overline{\beta + s} \setminus \beta$ with $k_1 < k_2 < \cdots < k_n$, to get an increasing chain of containment.

$$\beta \prec \overline{\beta + s}$$

$$\beta \cup \{k_1\} \prec \overline{\beta + s}$$

$$\beta \cup \{k_1, k_2\} \prec \overline{\beta + s}$$

$$\vdots$$

$$\beta \cup \{k_1, k_2, ..., k_n\} = \overline{\beta + s}$$

To motivate the following lemma, let us study the $(5,6)$-core bead diagram.

```
 4
 9    3
14    8    2
19   13    7    1
```

First, let us first consider all $(5,6)$-core beta-sets which have 3 elements: $\{1,2,3\}$,

$\{1,2,4\}$, $\{1,3,4\}$, $\{2,3,4\}$, $\{1,2,7\}$, $\{2,3,8\}$, and $\{3,4,9\}$. If we compare these beta-sets,

we find that all of them are contained in $\{3,4,9\}$ partition-wise. Similarly if we compare

all 4 element $(5,6)$-core beta-sets we find they are all contained partition-wise in

$\{9,4,3,2\}$, all 5 element $(5,6)$-core beta-sets are contained partition-wise in $\{9,8,4,3,2\}$, all

6 element $(5,6)$-core beta-sets are contained partition-wise in $\{14,9,8,4,3,2\}$, and so on.

If we arrange the maximal $n$ element $(5,6)$-core beta-sets in order they look like the following.

```
 4              4              4
9   3          9  [3]        [9] [3]
14  8  2       14  8   2      14  8   2
19 13 7 1      19 13  7 1     19 13  7 1

 4              4              4
[9] [3]        [9] [3]        [9] [3]
14  8 [2]       14 [8] [2]    [14] [8] [2]
19 13 7 1      19 13  7 1      19 13  7 1

 4              4              4
[9] [3]        [9] [3]        [9] [3]
[14] [8] [2]   [14] [8] [2]   [14] [8] [2]
19 13 7 [1]    19 13 [7] [1]  19 [13] [7] [1]

 4
[9] [3]
[14] [8] [2]
[19] [13] [7] [1]
```

These all take the same form: they look like a delta-set whose generator is in the left-most column, sometimes having elements missing from the bottom-most row. The following lemma shows that these maximal $n$ element beta-sets (called canonical forms) will always appear given the right conditions.

Here, $T_i$ denotes the $i$th triangular number; that is, $T_i = 1 + 2 + \cdots + i$.

## Lemma 4.4 (Canonical form of beta-sets)

Let $s,t$ be relatively prime positive integers, $n,\ i$ be positive integers, and $\beta$ be a beta-set, such that

($i$)     $i < \min\{s+1, \frac{n}{t}\}$ ;

($ii$)     $\beta \subset D := \{n - as - bt \mid a, b \geq 0; a + b \leq i\}$ ; and,

($iii$)     $\overline{\beta} \cap D = \beta$ (we will say $\beta$ is $(s,t)$-closed in $D$).

Now let $k$ be the index of the triangular number such that $T_k \leq |\beta| < T_{k+1}$, and $A$ be a subset of positive integers such that $A$ is empty if $T_k = |\beta|$, and $A = \{n - (i-k)s - bt \mid k + T_k - |\beta| < b \leq k\}$ otherwise. Then if we define a new beta-set $\gamma$ by $\gamma = \{n - (i - k + 1)s - as - bt \mid a, b \geq 0; a + b \leq k - 1\} \cup A$, we have that $\beta \prec \gamma \subset D$ and we call $\gamma$ the canonical form of $\beta$ in $D$.

## Proof of Lemma 4.4

First we observe that each element of $D$ can be written uniquely in the form $n - as - bt$. Suppose that $n - as - bt = n - a's - b't$, then $(a' - a)s = (b - b')t$. But $s, t$ are relatively prime, so $s \mid (b - b')$ and $t \mid (a' - a)$, but $a, a', b, b'$ are all non-negative integers smaller than $s$, which means that $a' - a = b - b' = 0$, and hence $n - as - bt$ is unique. Also note that the smallest element of $D$ is $n - it$, which, by the constraints on $i$, must be positive.

Hence the two constraints on $i$ imply that $D$ will not overlap itself in the extended bead diagram and that $D$ is fully contained (set-wise) in the extended bead diagram.

Now define $D_j = \{n - as - bt \mid a + b = j\}$. If we place $D$ in the extended $(s,t)$-core bead diagram and look at only one copy of $D$, then $D_j$ correspond to diagonals on the bead diagram intersected with $D$, with $D_0$ containing just $n$, $D_1$ containing the next two beads (the one above and to the right of $n$) and so on.

Since this is a very visual process, here are some examples from the (5,6)-case with $n = 19$, $i = 4$ (in this case, $D$ is the maximal (5,6)-core partition). Note that $D$ has the shape of an isosceles right triangle placed on one of its legs.

$$
\begin{array}{llll}
& 4 & 4 & 4 & \boxed{4} \\
& 9 \quad 3 & 9 \quad 3 & \boxed{9} \quad 3 & 9 \quad \boxed{3} \\
D_0: & 14 \quad 8 \quad 2 \quad\quad D_1: & \boxed{14} \quad 8 \quad 2 \quad\quad D_2: & 14 \quad \boxed{8} \quad 2 \quad\quad D_3: & 14 \quad 8 \quad \boxed{2} \\
& \boxed{19} \quad 13 \quad 7 \quad 1 & 19 \quad \boxed{13} \quad 7 \quad 1 & 19 \quad 13 \quad \boxed{7} \quad 1 & 19 \quad 13 \quad 7 \quad \boxed{1}
\end{array}
$$

Suppose $x = n - a_x s - b_x t \in \beta \cap D_j$ and $j < i$. Then $x - s = n - (a_x + 1)s - b_x t$ and $x - t = n - a_x s - (b_x + 1)t$ are elements of $D_{j+1}$ and since $j + 1 \le i$, $x - s, x - t \in D$. But that implies both $x - s, x - t$ are positive, since all elements of $D$ are positive. By the definition of $\overline{\beta}$, this implies that $x - s, x - t \in \overline{\beta}$. Since $\overline{\beta} \cap D = \beta$, $x - s, x - t$ are in $\beta$ as well.

Now suppose $\beta \cap D_j$ is non-empty, then for each $x \in \beta \cap D_j$, we know that $x - s \in \beta \cap D_{j+1}$. But if we let $x' = \min(\beta \cap D_j)$, then $x' - t \in \beta \cap D_{j+1}$, but $x' - t$ cannot be written as $x - s$ for any $x \in \beta \cap D_j$, since that would imply that $x' - t + s \in \beta \cap D_j$, which contradicts the fact that $x'$ is the minimal element of $\beta \cap D_j$. Therefore we know that $\beta \cap D_{j+1}$ must contain at least one more element than $\beta \cap D_j$ if $\beta \cap D_j$ is non-empty.

With this fact in mind, we create a canonical form of $\beta$ in two steps.

Step 1) We transform $\beta$ into $\beta'$ by sliding the elements of $\beta$ up along the diagonals. That is, if $\beta$ had 3 elements in $D_n$, then $\beta'$ has the 3 largest elements in $D_n$. We say then that $\beta'$ is diagonally top justified.

Suppose $x$ is the $m$th largest element of $j$th diagonal ($j < i$ again), then $x = n - (j + 1 - m)s - (m-1)t$. If $x \in \beta' \cap D_j$ then $\beta' \cap D_j$ must contain at least $m$ elements, so we know that $\beta' \cap D_{j+1}$ must contain at least $m+1$ elements. But the $m$th and $(m+1)$st largest elements of $\beta' \cap D_{j+1}$ are simply $n - (j + 2 - m)s - (m-1)t$ and $n - (j + 2 - (m+1))s - mt$, which are $x - s, x - t$ respectively. Thus $\beta'$ is $(s,t)$-closed in $D$.

Since either an element will stay the same or be moved to a higher value, we can apply Corollary 2.7 repeatedly to show that $\beta \prec \beta'$.

Staying with the example for $D$ above, suppose that

$$\beta = \begin{matrix} 4 & & & \\ 9 & \boxed{3} & & \\ 14 & 8 & \boxed{2} & \\ 19 & 13 & \boxed{7} & \boxed{1} \end{matrix} \quad .$$

Then,

$$\beta' = \begin{matrix} \boxed{4} & & & \\ \boxed{9} & \boxed{3} & & \\ 14 & 8 & \boxed{2} & \\ 19 & 13 & 7 & 1 \end{matrix} \quad .$$

But similarly, if

$$\beta = \begin{array}{ccccc} & \boxed{4} & & & \\ & 9 & 3 & & \\ & 14 & 8 & \boxed{2} & \\ & 19 & 13 & \boxed{7} & \boxed{1} \end{array}$$

Then again,

$$\beta' = \begin{array}{ccccc} & \boxed{4} & & & \\ & \boxed{9} & \boxed{3} & & \\ & 14 & 8 & \boxed{2} & \\ & 19 & 13 & 7 & 1 \end{array}$$

So note that the same $\beta'$ can come from two different choices for $\beta$.

Step 2) We transform $\beta'$ into $\beta''$ by a recursive method. We start with $j = 2$. For the remainder of the proof, let $x_j = x_j(\tilde{\beta}) = \min(\tilde{\beta} \cap D_j)$ where $\tilde{\beta}$ is whatever altered form of $\beta$ we are considering at the time. If $x_j + t \in \beta'$ then increase $j$ by 1 and start again. Otherwise consider $E(x_j) = \{x_j + t - s + bt \mid b > 0\} \cap (D \setminus \beta')$. If $E(x_j)$ is empty, increase $j$ by 1 and start again. If $E(x_j)$ contains at least one element, then remove $x_j$ from $\beta'$ and add $\min(E(x_j))$ in its place; slide the new element up along the diagonal as far as it will go; then reset $j$ to 2 and begin the second step again.

We finish when $j > i$.

This final beta-set we call $\beta''$.

We again want to show that $\beta''$ is $(s,t)$-closed in $D$. First we need to show that by removing $x_j$, we do not lose $(s,t)$-closure. But since we reached the point in Step 2 where we moved $x_j$, we must have that $x_j + t \notin \beta''$. Similarly, $x_j + s \notin \beta''$ because it is smaller than and in the same diagonal as $x_j + t$. Since neither of these elements are present, then removing $x_j$ will not affect whether or not $\beta''$ is $(s,t)$-closed in $D$.

Similarly we must also show that the new element we add in does not affect whether $\beta''$ is $(s,t)$-closed in $D$. We know that the element we add in is on the same diagonal as $\min(E(x_j))$. Since $x_j$ is in $\beta'$ (before the removal process) we know that $x_j + t - s$ is in $\beta'$ since $\beta'$ is diagonally top justified. But this means that $\min(E(x_j)) - t$ must be in $\beta'$ by the construction of $E(x_j)$. Now suppose when we shift $\min(E(x_j))$ along its diagonal, it moves up $m$ places, that is, we actually are replacing $x_j$ with $\min(E(x_j)) + m(t-s)$, then the diagonal containing $\min(E(x_j)) - t$ must also contain $\min(E(x_j)) + m(t-s) - t$ and $\min(E(x_j)) + (m+1)(t-s) - t = \min(E(x_j)) + m(t-s) - s$ because again, $\beta'$ is diagonally top justified. But that means that adding in $\min(E(x_j)) + m(t-s)$ does not affect whether $\beta''$ is $(s,t)$-closed in $D$, so in general step 2 preserves the fact that $\beta'$ is $(s,t)$-closed in $D$.

Also, $\min(E(x_j)) + m(t-s) > x_j$ so we can apply Corollary 2.7 as many times as we have to increase the value of an element to see that $\beta' \prec \beta''$.

For a brief example, if $\beta' = \begin{smallmatrix} & \boxed{4} & & \\ \boxed{9} & & \boxed{3} & \\ 14 & & 8 & & \boxed{2} \\ 19 & & 13 & & 7 & & \boxed{1} \end{smallmatrix}$ then we move step by step.

$$\beta' = \begin{array}{cccc} \boxed{4} & & & \\ \boxed{9} & \boxed{3} & & \\ 14 & 8 & \boxed{2} & \\ 19 & 13 & 7 & \boxed{1} \end{array} \quad => \quad \begin{array}{cccc} \boxed{4} & & & \\ \boxed{9} & \boxed{3} & & \\ 14 & 8 & \boxed{2} & \\ 19 & 13 & \boxed{7} & 1 \end{array} \quad => \quad \begin{array}{cccc} \boxed{4} & & & \\ \boxed{9} & \boxed{3} & & \\ 14 & \boxed{8} & \boxed{2} & \\ 19 & 13 & 7 & 1 \end{array} = \beta''$$

Now we want to know what kind of form $\beta''$ has. In order to escape the looping Step 2, $\beta''$ must satisfy two conditions: $\beta''$ must be diagonally top justified and for each diagonal $D_j$ either $x_j + t \in \beta''$ or $E(x_j) = \varnothing$. We also know that $\beta''$ must be $(s,t)$-closed in $D$.

Now consider $x_j \in \beta''$ where $x_j + t \notin \beta''$, then we know that $E(x_j) = \varnothing$. By the definition of $E(x_j)$, this implies that $\{x_j + t - s + bt \mid b > 0\} \cap D \subset \beta''$, but since $\beta''$ is $(s,t)$-closed, $x_j + t - s \in \beta''$ implies $x_j - s \in \beta''$, which implies $x_j - t - s \in \beta''$, and so on. Thus we have that $\{x_j - bt - s \mid b \in \mathbb{Z}\} \cap D \subset \beta''$. But again, since $\beta''$ is $(s,t)$-closed, $x_j - bt - s \in \beta''$ implies $x_j - bt - 2s \in \beta''$, which implies $x_j - bt - 3s \in \beta''$, and so on. Thus we have that $\{x_j - as - bt \mid a \in \mathbb{N}, b \in \mathbb{Z}\} \cap D \subset \beta''$.

We now show that $\beta''$ is as described in the statement of the Lemma. If $x_j + t \in \beta''$, then we must have that $x_{j-1} = x_j + t$, because if $\beta'' \cap D_j$ has $m$ elements, then $\beta'' \cap D_{j-1}$ can have at most $m-1$ elements, and $x_j + t$ is the $(m-1)$st element of $D_{j-1}$. If instead $x_j + t \notin \beta''$ but $x_j + t \in D$, then $x_{j-1} = x_j + 2t - s$, because $E(x_j) = \varnothing$ implies that $x_j + 2t - s \in \beta''$ and this time $\beta'' \cap D_{j-1}$ can have at most $m-2$ elements, but again $x_j + 2t - s$ is the $(m-2)$nd element of $D_{j-1}$ and so $D_{j-1}$ has exactly $m-2$ elements.

Given this, let $y = n - a_y s - b_y t$, where $a_y = \min\{a_j \mid \hat{j} < j \le i\}$ and $b_y = \min\{b_j \mid \hat{j} < j < i; a_j = a_y\}$, and $\hat{j}$ represents the smallest index for which $x_j$ does not exist (because the intersection is empty). Then, $y \in \beta''$, but more than that, $y$ characterizes $\beta''$. To see this, we know for some $l$, $y = x_l$. By the previous paragraph, $x_j = y - (j - l)t$ for $j > l$ and $x_j = y - (j - l + 1)t - s$ for $\hat{j} < j < l$.

Thus, $\beta'' = \{n - as - bt \mid a > a_y, a + b \le i\} \cup \{y, y - t, y - 2t, ..., y - (i - a_y - b_y)t\}$ and since $|\beta| = |\beta''|$ and since for different values of $y$, $\beta''$ will contain a different number of elements, we can see that $\beta''$ is dependent only upon the number of elements in the original beta-set $\beta$. We can thus rewrite $\beta''$ as $\{n - (a_y + 1)s - as - bt \mid a + b \le i - a_y - 1\} \cup A$, where the first set has $T_k$ $(k = i - a_y - 1)$ elements and the second set is the same set $A$ we refer to in the statement of Lemma 4.4. But this means $\gamma = \beta''$ and hence $\beta \prec \gamma$.

Q.E.D.

We will now prove the maximal theorem in the case when $t = s + 1$.

**Proof of the $(s, s+1)$ case**

Let $\beta \subset \beta_{s,s+1}$ be a beta-set. Then note that $\beta$ and $D = \beta_{s,s+1}$ satisfy the conditions for the previous lemma. Therefore we know that $\beta \prec \gamma$, where $\gamma$ has the form $\gamma = \{s(s+1) - s - (s+1) - (i - k + 1)s - as - b(s+1) \mid a, b \ge 0; a + b \le k - 1\} \cup A$. If $A$ is empty, we say $\gamma$ is a type I canonical beta-set and we say $\gamma$ is type II otherwise.

Now we will show all canonical beta-sets are contained partition-wise in $\beta_{s,s+1}$. First, note that all type I canonical beta-sets can be rewritten as $\overline{\{ns-1\}}$. By Theorem 3.1, $\overline{\{ns-1\}} \prec \overline{\{ns-1\}+s} = \overline{\{(n+1)s-1\}}$. So each type I canonical beta-set is contained partition-wise in the next largest type I canonical beta-set. Since $\beta_{s,s+1}$ is itself a type I canonical beta-set, all type I canonical beta-sets are contained partition-wise in $\beta_{s,s+1}$.

All type II canonical beta-sets are of the form $\overline{\{ns-1,k\}}$. These can be rewritten as $\overline{\{ns-1,k\}} = \overline{\{ns-1\}} \cup \{k, k-(s+1), k-2(s+1),...,k-j(s+1)\}$ where $j$ is as large as it can be without $k-j(s+1)$ being negative. Then as we discussed earlier, we simply apply Lemma 4.3 to $\overline{\{ns-1\}}$ and $\overline{\{(n+1)s-1\}}$. This gives us $\overline{\{ns-1\}} \cup \{k-j(s-1)\} \prec \overline{\{(n+1)s-1\}}$. Applying Lemma 4.3 to $\overline{\{ns-1\}} \cup \{k-j(s-1)\}$ and $\overline{\{(n+1)s-1\}}$ gives $\overline{\{ns-1\}} \cup \{k-(j-1)(s-1),k-j(s-1)\} \prec \overline{\{(n+1)s-1\}}$. And so we can repeat Lemma 4.3 $j$ times to find that $\overline{\{ns-1,k\}} = \overline{\{ns-1\}} \cup \{k,k-(s+1),k-2(s-1),...,k-j(s-1)\} \prec \overline{\{(n+1)s-1\}}$. Thus every type II canonical beta-set is contained partition-wise in the next largest (under set containment) type I canonical beta-set which we've already shown is contained partition-wise in $\beta_{s,s+1}$.

Thus, all $(s,s+1)$-core beta-sets are contained partition-wise in their respective $\beta''$, and each $\beta''$ is contained partition-wise in $\beta_{s,s+1}$, so the theorem holds.

Q.E.D.

We can extend this proof to cover the $(s, ks+1)$-core case with only minimal changes to the proof. We would need to have a new notion of a canonical beta-set and its diagonals, as well as a slightly different method of constructing the canonical form.

The diagonals we can redefine as $D_j = \{n - as - bt \mid a + kb = j\}$ and the second stage of creating the canonical beta-set we define by pushing an element from one diagonal to a smaller indexed diagonal provided the beta-set remains $(s, sk+1)$-closed in $D$.

However, we cannot use this method to prove the maximal theorem for generic $t$, as the diagonals as we originally defined them no longer increase in size as we increase the index, and there is no way of defining them to do so in a regular fashion.

## 5. Proving the general case

We require several more lemmas to be able to prove the general case of the maximal theorem. The first of these provides sufficient conditions for preserving partition-wise containment when adding the same constant to both beta-sets.

**Lemma 5.1** Given two beta-sets $\beta, \gamma$ such that $\beta \prec \gamma$, a positive integer $k$, and two sets $A, B \subset \{1, 2, ..., k\}$ such that $|A| \geq |B|$, and $|\gamma| \geq |(\beta + k) \cup A|$, we then have also that $(\beta + k) \cup A \prec (\gamma + k) \cup B$.

**Proof of Lemma 5.1**

To show $(\beta + k) \cup A \prec (\gamma + k) \cup B$, we will use comparison by separation (Proposition 2.8).

We want to show partition-wise containment of the first $|\beta|$ elements; however, this is simple. For $1 \leq i \leq |\beta|$,

$$P((\beta + k) \cup A)_i$$

$$= P(\beta + k)_i - |A|$$

$$= P(\beta)_i + k - |A|$$

$$\leq P(\beta)_i + k - |B| \qquad \text{because } |A| \geq |B|$$

$$\leq P(\gamma)_i + k - |B| \qquad \text{because } \beta \prec \gamma$$

$$= P(\gamma + k)_i - |B|$$

$$= P((\gamma + k) \cup B)_i.$$

The remaining elements to be compared are those in $A$; however, because $|\gamma| \geq |(\beta + k) \cup A|$, we know that the elements of $A$ are being compared with elements in $(\gamma + k)$, not with the elements in $B$.

Now we can simply take worst case scenario and assume the elements of $A$ are as large as possible and the elements of $(\gamma + k)$ are as small as possible. Let $(\gamma + k)_i^j$ denote the subset of $(\gamma + k)$ containing the $i$th through $j$th elements. Then in the worst case, $A = \{k, k-1, ..., k - |A| + 1\}$ and $(\gamma + k)_{|\beta|+1}^{|\gamma|} = (\gamma)_{|\beta|+1}^{|\gamma|} + k = \{|\gamma| - |\beta|, |\gamma| - |\beta| - 1, ..., 1\} + k = \{|\gamma| - |\beta| + k, |\gamma| - |\beta| - 1 + k, ..., 1 + k\}$. However, $(\gamma + k)$ is followed by the set $B$ and so the partition is reduced by $|B|$. Thus we are comparing $P(A) = P(\{k, k-1, ..., k - |A| + 1\}) = (k - |A| + 1, k - |A| + 1, ..., k - |A| + 1)$ with $P(\gamma + k)_{|\beta|+1}^{|\gamma|} - |B| = (1 + k - |B|, 1 + k - |B|, ..., 1 + k - |B|)$. We then have partition-wise containment provided that $1 + k - |B| \geq k - |A| + 1$, but this is true since $|A| - |B| \geq 0$.

Thus we have that $(\beta + k) \cup A \prec (\beta + k) \cup B$ by Proposition 2.8.

Q.E.D.

We will often use Lemma 5.1 in the following form.

**Corollary to Lemma 5.1**  Given two $(s,t)$-closed beta-sets $\beta, \gamma$ and a positive integer $k$, such that $\beta \prec \gamma$, $\left|\overline{\beta+k}\right| - |\beta| \geq \left|\overline{\gamma+k}\right| - |\gamma|$, and $|\gamma| \geq \left|\overline{\beta+k}\right|$, we then have that $\overline{\beta+k} \prec \overline{\gamma+k}$ as well.

**Proof of Corollary to Lemma 5.1**

This follows from Lemma 5.1 if we let $\overline{\beta+k} = (\beta+k) \cup A$  and  $\overline{\gamma+k} = (\gamma+k) \cup B$.

Q.E.D.

For our next lemma we need to have some way of glossing over that $+s$ term in the second generalization.  The following lemma says in essence that when we get close enough for the $+s$ error term to bother us we already have partition-wise containment by a different method.

**Lemma 5.2**  For any $(s,t)$-closed delta-set $\Delta$ and any $(s,t)$-closed beta-set $\beta$, such that $\beta \subset \Delta$ and $|\Delta| - |\beta| \leq w(\Delta)$, we have that $\beta \prec \Delta$.

**Proof of Lemma 5.2**

Let $D$ be the subset of $\Delta$ defined by the following rule: if $\Delta = \overline{\{n\}}$, then $D = \{n - as - bt \mid a, b \geq 0, a + b < w(\Delta)\}$.  Thus $D$ has the shape of an isosceles right triangle lying on one of its legs.

Here's an example of $D$ for (8,11) and $\Delta$ being $\beta_{8,11}$ (elements of $D$ are blocked in). Note that $D$ is seven elements wide and seven elements tall.

```
 5
13     2
21    10
29    18     7
37    26    15     4
45    34    23    12     1
53    42    31    20     9
61    50    39    28    17     6
69    58    47    36    25    14     3
```

All of the elements in $\Delta \setminus \beta$ must be elements of $D$. This is true because we can again think of $D$ as a union of diagonals $D_j = \{n - as - bt \mid a,b \geq 0, a + b = j\}$ where $j$ ranges from 0 to $w(\Delta) - 1$. Since $\beta$ is $(s,t)$-closed, we have the rule that if $x \in \beta$ implies $x - s, x - t \in \beta$ (if the values are positive). But we also have a reversal of this rule by Proposition 2.9: if $x \notin \beta$ then $x + s, x + t \notin \beta$. Thus if $x \in D_j$ but $x \notin \beta$ then at least one of $x + s, x + t$ is in $D_{j-1}$ and that same element cannot be in $\beta$. But by the same argument, there must be an element in $D_{j-2}$ that isn't in $\beta$, and so on.

So now consider $x \in \Delta \setminus (D \cup \beta)$, then there exist $a,b$ such that $x = n - as - bt$ (there might exist more than one such choice of $a$ and $b$ but that will not matter). Then consider the set $\{n - as - bt, n - (a-1)s - bt, ..., n - s - bt, n - bt, n - (b-1)t, ..., n\}$. By Proposition 2.9, this set is disjoint from $\beta$, but at the same time it is contained set-wise in $\Delta$. But $\{n-as-bt, n-(a-1)s-bt, ..., n-s-bt, n-bt, n-(b-1)t, ..., n\}$ contains at least $w(\Delta)+1$ elements since $x \notin D$ implies that $a+b > w(\Delta)$.

To follow our example above, suppose 31 were not an element of $\beta$, then neither would 39, 47, 58, or 69 (among others) be and each of these is in a diagonal one index smaller than the last. Similarly if 1 were not an element of $\beta$, then neither would be 9, 17, 25, 36, 47, 58, or 69 (among others), giving us at least 8 elements in $|\Delta| - |\beta|$ whereas $w(\Delta) = 7$.

Then we note that $D \cap \beta$ and $D$ satisfy the conditions for Lemma 4.4. Thus we have a canonical form of $D \cap \beta$ in $D$, which we call $\gamma$. Then we can see that $\beta = (D \cap \beta) \cup (\Delta \setminus D) \prec \gamma \cup (\Delta \setminus D)$. This last beta-set is $(s,t)$-closed because the $(s,t)$-closure of $\gamma$ is $\gamma \cup (\Delta \setminus D)$.

Let us demonstrate this with a full example.

We start with a beta-set $\beta$ as a part of our delta-set $\beta_{8,11}$.

$$\beta = \begin{array}{ccccccc}
\boxed{5} & & & & & & \\
\boxed{13} & \boxed{2} & & & & & \\
\boxed{21} & \boxed{10} & & & & & \\
\boxed{29} & \boxed{18} & \boxed{7} & & & & \\
\boxed{37} & \boxed{26} & \boxed{15} & \boxed{4} & & & \\
45 & \boxed{34} & \boxed{23} & \boxed{12} & \boxed{1} & & \\
53 & \boxed{42} & \boxed{31} & \boxed{20} & \boxed{9} & & \\
61 & 50 & \boxed{39} & \boxed{28} & \boxed{17} & \boxed{6} & \\
69 & 58 & 47 & \boxed{36} & \boxed{25} & \boxed{14} & \boxed{3}
\end{array}$$

Then we consider $D \cap \beta =$

```
[21]
[29] [18]
[37] [26] [15]
 45  [34] [23] [12]
 53  [42] [31] [20] [9]
 61   50  [39] [28] [17] [6]
 69   58   47  [36] [25] [14] [3]
```

Now we apply the canonical form creation process to create $\gamma$.

```
[21]
[29] [18]
[37] [26] [15]
 45  [34] [23] [12]
 53  [42] [31] [20] [9]
 61   50  [39] [28] [17] [6]
 69   58   47  [36] [25] [14] [3]
```

```
[21]
[29] [18]
[37] [26] [15]
[45] [34] [23] [12]
 53  [42] [31] [20] [9]
 61   50  [39] [28] [17] [6]
 69   58   47   36  [25] [14] [3]
```

 

**Diagram 1**

| [21] | | | | | | |
|---|---|---|---|---|---|---|
| [29] | [18] | | | | | |
| [37] | [26] | [15] | | | | |
| [45] | [34] | [23] | [12] | | | |
| [53] | [42] | [31] | [20] | [9] | | |
| 61 | 50 | 39 | [28] | [17] | [6] | |
| 69 | 58 | 47 | 36 | [25] | [14] | [3] |

**Diagram 2**

| [21] | | | | | | |
|---|---|---|---|---|---|---|
| [29] | [18] | | | | | |
| [37] | [26] | [15] | | | | |
| [45] | [34] | [23] | [12] | | | |
| [53] | [42] | [31] | [20] | [9] | | |
| 61 | 50 | [39] | [28] | [17] | [6] | |
| 69 | 58 | 47 | 36 | 25 | [14] | [3] |

**Diagram 3**

| [21] | | | | | | |
|---|---|---|---|---|---|---|
| [29] | [18] | | | | | |
| [37] | [26] | [15] | | | | |
| [45] | [34] | [23] | [12] | | | |
| [53] | [42] | [31] | [20] | [9] | | |
| 61 | [50] | [39] | [28] | [17] | [6] | |
| 69 | 58 | 47 | 36 | 25 | 14 | [3] |

**Diagram 4**  = γ

| [21] | | | | | | |
|---|---|---|---|---|---|---|
| [29] | [18] | | | | | |
| [37] | [26] | [15] | | | | |
| [45] | [34] | [23] | [12] | | | |
| [53] | [42] | [31] | [20] | [9] | | |
| [61] | [50] | [39] | [28] | [17] | [6] | |
| 69 | 58 | 47 | 36 | 25 | 14 | 3 |

Now we reattach all the other elements to get the following.

$$\gamma \cup (\Delta \setminus D) =$$

|   5 |     |     |     |     |    |   |
|-----|-----|-----|-----|-----|----|---|
|  13 |   2 |     |     |     |    |   |
|  21 |  10 |     |     |     |    |   |
|  29 |  18 |   7 |     |     |    |   |
|  37 |  26 |  15 |   4 |     |    |   |
|  45 |  34 |  23 |  12 |   1 |    |   |
|  53 |  42 |  31 |  20 |   9 |    |   |
|  61 |  50 |  39 |  28 |  17 |  6 |   |
|  69 |  58 |  47 |  36 |  25 | 14 | 3 |

But this construction implies that $\gamma \cup (\Delta \setminus D)$ consists of all elements of $\Delta$ except for some or all elements in the bottom row. Again, as in the proof of the $(s,s+1)$ case, this means that we have $\gamma \cup (\Delta \setminus D) \prec \Delta$ and hence $\beta \prec \Delta$.

Q.E.D.

Since we have this lemma already at our disposal we will show now how the second generalization implies the first generalization.

**Theorem 3.2 (First generalization of the maximal theorem)**

Let $\Delta$ be any $(s,t)$-closed delta-set and $\beta$ be any $(s,t)$-closed beta-set such that $\beta \subset \Delta$. Then $\beta \prec \Delta$.

**Proof of the first generalization of the maximal theorem**

Given an $(s,t)$-closed delta-set $\Delta$, consider a beta-set $\beta \subset \Delta$.

First case: $|\Delta| - |\beta| > w(\Delta)$

Since $|\Delta - s| = |\Delta| - w(\Delta)$, we have that $|\Delta - s| > |\beta| = \left| M\Delta_{|\beta|} \right|$, but $\Delta - s$ and $M\Delta_{|\beta|}$

are delta-sets so this implies that the generator of $\Delta - s$ is larger than the generator of

$M\Delta_{|\beta|}$ (see Theorem 3.1.$vi$). Thus the generator of $\Delta$ is larger than the generator of

$\overline{M\Delta_{|\beta|} + s}$. But that again implies that $\overline{M\Delta_{|\beta|} + s} \prec \Delta$.

By the second generalization, we have $\beta \prec \overline{M\Delta_{|\beta|} + s} \prec \Delta$.

Second case: $|\Delta| - |\beta| \leq w(\Delta)$

In this case, just apply Lemma 5.2 to see that $\beta \prec \Delta$.

Q.E.D.

Before we prove the second generalization we wish to prove some lemmas to simplify some steps in the proof of the second generalization.

**Lemma 5.3** If $a < b$ then $w(\overline{\{a\}}) \leq w(\overline{\{b\}})$ and $h(\overline{\{a\}}) \leq h(\overline{\{b\}})$. Furthermore, given a

$w \leq s$, there exists a delta-set $\Delta$ with $w(\Delta) = w$. Alternately, given a $h \leq t$, there exists

a delta-set $\Delta$ with $h(\Delta) = h$.

**Proof of Lemma 5.3**

Note first that $\overline{\{b\}} = \overline{\overline{\{a\}} + (b-a)} = \left(\overline{\{a\}} + (b-a)\right) \cup A$, where $A \subset \{1, 2, ..., b-a\}$ as we know from Lemma 4.1.

Then suppose we reduce the elements of $\overline{\{b\}}$ modulo $s$ to find the width. This is equivalent then to reducing the elements of $\left(\overline{\{a\}} + (b-a)\right) \cup A$ modulo $s$, but $\overline{\{a\}} + (b-a)$ has the same number of residue classes as $\overline{\{a\}}$. In fact if $c$ is a residue class for $\overline{\{a\}}$, then $c + b - a$ is a residue class for $\overline{\{a\}} + (b-a)$. Therefore $\overline{\{b\}}$ must contain at least as many residue classes modulo $s$ as $\overline{\{a\}}$, therefore $w(\overline{\{a\}}) \leq w(\overline{\{b\}})$.

By the same argument now modulo $t$, we have that $h(\overline{\{a\}}) \leq h(\overline{\{b\}})$.

To prove that there exist delta-sets with any possible width (height), it is enough to note that by Theorem 3.1 and the work above, $0 \leq w(\overline{\{a+1\}}) - w(\overline{\{a\}}) \leq 1$, and that there exist delta-sets with width equal to 0 and delta-sets with width equal to $s$.

Q.E.D.

**Lemma 5.4** Let $\Delta^1, \Delta^2, ..., \Delta^k$ be a set of disjoint $(s,t)$-closed delta-sets such that $k > 2$ and at least one $\Delta^i$ has width greater than 1. Let $n = \left|\Delta^1\right| + \left|\Delta^2\right| + \cdots + \left|\Delta^k\right|$. Then

$h(M\Delta_n) < h(\Delta^1) + h(\Delta^2) + \cdots + h(\Delta^k)$.

**Proof of Lemma 5.4**

Let $\Gamma$ be the smallest delta-set with respect to partition containment whose height equals $h(\Delta^1) + h(\Delta^2) + \cdots + h(\Delta^k)$. Then following the Euclidean algorithm, there exist positive integers $q, r$ such that $0 \leq r < s$ and $t = qs + r$.

Suppose we have an arbitrary delta-set with generator small enough that each element has a unique representation of the form $n - as - bt$, where $n$ is the generator, and we wish to know how many elements are in two adjoining columns. By adjoining columns we mean that the largest element in the column on the left should be $t$ larger than the largest element in the column on the right, like so.

$$
\begin{array}{ll}
a & \\
\vdots & \\
a + (m_a - m_b - 1)s & \\
a + (m_a - m_b)s & b \\
\vdots & \vdots \\
a + (m_a - 1)s & b + (m_b - 1)s \\
a + m_a s & b + m_b s
\end{array}
$$

In this diagram, $b + m_b s = a + m_a s - t$, and $0 < a, b \leq s$. The column on the left has $m_a + 1$ elements and the column on the right has $m_b + 1$ elements.

So suppose $a \leq r$, then $a + qs \leq r + qs = t$, thus $a + (q+1)s - t$ is an element of the right column that is smaller than $s$. So $b = a + (q+1)s - t$. In this case, $m_b = m_a - (q+1)$, so the left column has $q + 1$ more elements.

Suppose $a > r$, then $a + qs > r + qs = t$, so $a + qs - t$ is an element of the right column and it is also smaller than $s$. So $b = a + qs - t$. In this case, $m_b = m_a - q$, so the left column has $q$ more elements.

Thus the difference in size between two adjoining columns is at most $q+1$ and at least $q$.

So returning to $\Gamma$ and $\Delta^1, \Delta^2, ..., \Delta^k$, let us compare how many elements are in the $i$th column from the left in $\Gamma$ and how many elements are in each of the $i$th columns from the right in all the $\Delta^1, \Delta^2, ..., \Delta^k$.

The 0th column in $\Gamma$ has exactly $h(\Gamma) = h(\Delta^1) + h(\Delta^2) + \cdots + h(\Delta^k)$ elements. Similarly the 0th columns in the $\Delta^1, \Delta^2, ..., \Delta^k$ have $h(\Delta^1) + h(\Delta^2) + \cdots + h(\Delta^k)$.

Then the $i$th column of $\Gamma$ has at least $h(\Gamma) - (q+1)i$ elements and the $i$th columns of $\Delta^1, \Delta^2, ..., \Delta^k$ have at most $(h(\Delta^1) - qi) + (h(\Delta^2) - qi) + \cdots + (h(\Delta^k) - qi) = h(\Gamma) - kqi$ elements. But since $k > 2$, $h(\Gamma) - kqi < h(\Gamma) - (q+1)i$.

Since $\Gamma$ has more elements than $\Delta^1, \Delta^2, ..., \Delta^k$ in each column after the 0th, we have that $|\Gamma| > |\Delta^1| + |\Delta^2| + \cdots + |\Delta^k| = n$. Thus by Theorem 3.1, the generator of $\Gamma$ must be larger than the generator of $M\Delta_n$. But $\Gamma$ was the smallest delta-set (with respect to partition containment) with height equal to $h(\Delta^1) + h(\Delta^2) + \cdots + h(\Delta^k)$, and since $M\Delta_n$ is smaller still, we must have that $h(M\Delta_n) < h(\Delta^1) + h(\Delta^2) + \cdots + h(\Delta^k)$.

Q.E.D.

**Lemma 5.5** Let $\beta$ be an $(s,t)$-closed beta-set containing $n$ elements such that $w(\beta) = w(M\Delta_n) < s$ and $h(\beta) = h(M\Delta_n)$. If we define $D$ as the smallest delta-set with respect to partition-containment which has $\beta$ as a subset, then $w(D) = w(\beta)$ and $h(D) = h(\beta)$.

**Proof of Lemma 5.5**

Suppose that the height and width of $D$ and $\beta$ are not the same. If $D$ has a residue class $a$ modulo $s$ not in $\beta$, where $1 \leq a \leq s$ then $\beta$ must also be missing the residue class $a$ modulo $t$ since it is $(s,t)$-closed. Therefore it is sufficient to show that the heights of $D$ and $\beta$ are equal.

For every residue class modulo $t$ present in $D$ but not $\beta$, we know that the representative $j$ ($1 \leq j \leq t$) of that residue class must be an element of $D \backslash \beta$. If $D = \overline{\{d\}}$ then for some $a_j, b_j$, $j = d - a_j s - b_j t$. But since $j \notin \beta$ no element of the form $j + as + bt$ can be an element of $\beta$ either. Thus, $\beta \subset D \backslash \{d - as - bt \mid 0 \leq a \leq a_j; 0 \leq b \leq b_j\}$

$= \overline{\{d - (a_j + 1)s\}} \cup \overline{\{d - (b_j + 1)t\}}$.

These last two sets are disjoint, because if $x \in \overline{\{d - (a_j + 1)s\}} \cap \overline{\{d - (b_j + 1)t\}}$, then on the one hand $x$ must be of the form $d - (a_j + 1 + a)s - bt$ in order to be in the first set, and on the other hand $x$ must also be of the form $d - (b_j + 1 + b)t - as$ to be in the second set, thus $x = d - (a_j + 1 + a)s - (b_j + 1 + b)t$. But remember that $j = d - a_j s - b_j t \leq t$ so $x \leq 0$ which is impossible. Thus the two sets are disjoint.

If we find that $h(\overline{\{d - (a_j + 1)s\}} \cup \overline{\{d - (b_j + 1)t\}}) > h(\beta)$ still, then we can repeat the process of the previous two paragraphs to find three disjoint delta-sets which contain $\beta$ set-wise. Repeating this process enough times, we will eventually find a disjoint set of delta-sets $\Delta^1, \Delta^2, ..., \Delta^k$ such that $\beta \subset \bigcup \Delta^k$, $h(\beta) = h(\Delta^1) + h(\Delta^2) + \cdots + h(\Delta^k)$, and $|\beta| \leq |\Delta^1| + |\Delta^2| + \cdots + |\Delta^k|$.

Suppose that $\beta$ and the $\Delta^1, \Delta^2, ..., \Delta^k$ satisfy the conditions of the previous lemma. If we let $n' = |\Delta^1| + |\Delta^2| + \cdots + |\Delta^k|$, then by the previous lemma $h(M\Delta_{n'}) < h(\Delta^1) + h(\Delta^2) + \cdots + h(\Delta^k)$, but at the same time, since $n \leq n'$ we have that $h(M\Delta_n) \leq h(M\Delta_{n'})$. Putting this all together gives $h(M\Delta_n) < h(\beta)$, which contradicts our original assumptions.

Suppose then that each of our $\Delta^1, \Delta^2, ..., \Delta^k$ has width 1. Then $\beta = \bigcup \Delta^k$, $h(\beta) = n$, and $w(\beta) = k$. If $M\Delta_n$ has width 1 as well, then we must have that $k = 1$ and hence $\beta = M\Delta_n - c$ for some integer $c$. If $M\Delta_n$ has width greater than 1, then $h(M\Delta_n) < n$, which contradicts our assumptions.

Finally, suppose then that $k \leq 2$ and the delta-sets have arbitrary width. If, in fact, $k = 1$, then the theorem is trivial, so let us assume $k = 2$. Here we just modify the argument of the previous lemma. If we let $\Gamma$ be the smallest delta-set whose height is at least $h(\Delta_1) + h(\Delta_2)$ and whose width is at least $w(\Delta_1) + w(\Delta_2)$, then in this case we have $h(\Gamma) - 2qi \leq h(\Gamma) - (q+1)i$ so just by comparing columns we cannot tell that $\Gamma$ has more than $n$ elements as we did in the previous lemma. But since here $w(\Gamma) \geq w(\Delta^1) + w(\Delta^2)$, we have that in the $(\max(w(\Delta^1), w(\Delta^2)) + 1)$st column $\Gamma$ must have at least one element while $\Delta^1, \Delta^2$ do not. Hence $|\Gamma| > |\Delta^1| + |\Delta^2| \geq |\beta|$ so we again have a contradiction.

Q.E.D.

So now we will prove the second generalization, and in so doing prove the maximal theorem itself.

**Theorem 4.2 (Second generalization of the maximal theorem)**

Let $M\Delta_n$ be the largest (under partition containment) $(s,t)$-closed delta-set containing $n$ elements. Then for any $(s,t)$-closed beta-set $\beta$ containing $n$ elements,

$$\beta \prec \overline{M\Delta_n + s} \,.$$

**Proof of second generalization of the maximal theorem**

We proceed by induction on $n$ and start by noticing that the generalization is trivially true if $\beta$ is a delta-set, since in this case $\beta \prec M\Delta_{|\beta|} \prec \overline{M\Delta_{|\beta|} + s}$ (see Theorem 3.1.$v$).

If $n = 1$, then note that all beta-sets containing just one element are, in fact, delta-sets, so the theorem holds for this case.

Now assume the generalization is true up to a given $n$ and we will prove it true for $n+1$.

Let $\beta$ be a beta-set of $n+1$ elements. If $\beta$ is a delta-set then we already know that the second generalization holds, so without loss of generality we will assume $\beta$ is not a delta-set. Then we break the proof into several cases based on the relative structures of $\beta$ and $M\Delta_{n+1}$.

Before we analyze cases specifically we must explain why some cases do not appear. Notably it is impossible for $w(\beta) < w(M\Delta_{n+1})$ or $h(\beta) < h(M\Delta_{n+1})$.

Suppose $w(\beta) < w(M\Delta_{n+1})$. If $w(M\Delta_{n+1}) = 1$, then clearly this is impossible so let $w(M\Delta_{n+1}) > 1$.

By the logic of the previous lemma's proof, we can find a set of disjoint delta-sets $\Delta^1, \Delta^2, ..., \Delta^k$ such that $\beta \subset \Delta^1 \cup \Delta^2 \cup \cdots \cup \Delta^k$ and $w(\beta) = w(\Delta^1) + w(\Delta^2) + \cdots + w(\Delta^k)$. If we want to explicitly calculate these delta-sets, then let $\Delta$ be the smallest delta-set with respect to partition containment that has $\beta$ as a subset, and let $g$ be its generator. Then for each representative for a given column residue class modulo $s$, $j = g - a_j s - b_j t \in \{1, 2, ..., s\}$ such that $j \in \Delta \backslash \beta$, remove the set $\{g - as - bt \mid 0 \le a \le a_j; 0 \le b \le b_j\}$ from $\Delta$. What we are left with is our set $\Delta^1, \Delta^2, ..., \Delta^k$. The order of these delta-sets does not matter.

Suppose initially that $k > 1$.

Now consider the generator of $\Delta^i$, and call it $g_i$. Let $\Gamma = \overline{\{g_k + (w(\beta) - w(\Delta^k))t\}}$. By construction $w(\Gamma) = w(\beta)$. But we also have that $|\Gamma| \ge |\Delta^1| + |\Delta^2| + \cdots + |\Delta^k|$ since for any $g_i - as - bt \in \Delta^i$, there is a corresponding element in $\Gamma$, namely the element $g_k + (w(\Delta^i) + w(\Delta^{i+1}) + \cdots + w(\Delta^{k-1}))t - as - bt$. This element is positive since for any $i \ne k$, $g_k + (w(\Delta^i) + w(\Delta^{i+1}) + \cdots + w(\Delta^{k-1}))t > w(\Delta^i)t \ge g_i$.

These elements are in fact unique corresponding elements. To see this, first note that since each delta-set $\Delta^i$ has width strictly smaller than $s$ we must have unique representation of its elements. If we have distinct elements $g_i - a_1 s - b_1 t, g_i - a_2 s - b_2 t \in \Delta^i$, then $g_k + (w(\Delta^i) + w(\Delta^{i+1}) + \cdots + w(\Delta^{k-1}))t - a_1 s - b_1 t \ne g_k + (w(\Delta^i) + w(\Delta^{i+1}) + \cdots + w(\Delta^{k-1}))t - a_2 s - b_2 t$ because $w(\Gamma) < s$ implies that each element has a unique representation in $\Gamma$ as well. Similarly if $i < j$ and $g_i - a_1 s - b_1 t \in \Delta^i, g_j - a_2 s - b_2 t \in \Delta^j$, then first note that $w(\Delta^{i+1}) + \cdots + w(\Delta^{k-1}) < w(\Delta^i) + w(\Delta^{i+1}) + \cdots + w(\Delta^{k-1}) - b_1 \leq w(\Delta^i) +$

$w(\Delta^{i+1}) + \cdots + w(\Delta^{k-1})$ and $w(\Delta^{j+1}) + \cdots + w(\Delta^{k-1}) < w(\Delta^j) + w(\Delta^{j+1}) + \cdots + w(\Delta^{k-1}) - b_1 \leq$

$w(\Delta^j) + w(\Delta^{j+1}) + \cdots + w(\Delta^{k-1})$. Thus we have that $(w(\Delta^i) + \cdots + w(\Delta^{j-1}) - b_1 +$

$b_2)t \neq (a_1 - a_2)s$ since $(w(\Delta^i) - b_1) + \cdots + w(\Delta^{j-1}) + b_2 < w(\Delta^i) + \cdots + w(\Delta^j) \leq w(\beta) < s$ and

$s, t$ are relatively prime, so $s$ cannot divide $(w(\Delta^i) + \cdots + w(\Delta^{j-1}) - b_1 + b_2)t$. Then $g_k +$

$(w(\Delta^i) + w(\Delta^{i+1}) + \cdots + w(\Delta^{k-1}))t - a_1 s - b_1 t \neq g_k + (w(\Delta^j) + w(\Delta^{j+1}) + \cdots + w(\Delta^{k-1}))t -$

$a_2 s - b_2 t$ since the coefficients of $t$ are distinct and $\Gamma$ has unique representation.

But now $w(\Gamma) < w(M\Delta_{n+1})$, which implies that the generator of $\Gamma$ (call it $g_\Gamma$) is smaller than the generator of $M\Delta_{n+1}$ (call it $g_{M\Delta_{n+1}}$) since both are delta-sets. Hence we know that for any element $g_\Gamma - as - bt \in \Gamma$, there is a corresponding element $g_{M\Delta_{n+1}} -$

$as - bt \in M\Delta_{n+1}$, but at the same time, the element $g_{M\Delta_{n+1}} - (w(M\Delta_{n+1}) - 1)t \in M\Delta_{n+1}$ does

not have a corresponding element in $\Gamma$ since $g_\Gamma - (w(M\Delta_{n+1}) - 1)t < 0$. Thus $|\Gamma| <$

$|M\Delta_{n+1}|$ and hence $|\beta| \leq |\Delta^1| + |\Delta^2| + \cdots + |\Delta^k| = |\Gamma| < |M\Delta_{n+1}|$, which is a contradiction.

Similarly if $k = 1$, then since $\beta$ is not a delta-set by assumption we have that

$|\beta| < |\Delta^1| < |M\Delta_{n+1}|$.

Thus we must have that $w(\beta) \geq w(M\Delta_{n+1})$.

The argument for the heights is proved similarly.

So now we move on to the individual cases.

First case: $w(\beta) > w(M\Delta_{n+1})$.

Instead of doing an immediate comparison, push $\beta$ and $M\Delta_{n+1}$ upwards and consider $\beta - s$ and $M\Delta_{n+1} - s$. Since a reduction by $s$ pushes the set upwards, all elements smaller than or equal to $s$ will disappear, and so this action shrinks the set by a number of elements equal to its width. But we have assumed that sets with fewer than $n+1$ elements satisfy the generalization, so we know that $\beta - s \prec \overline{M\Delta_{n+1-w(\beta)} + s}$, but at the same time $\left| M\Delta_{n+1-w(\beta)} \right| = |\beta - s| = |\beta| - w(\beta) < |M\Delta_{n+1}| - w(M\Delta_{n+1}) = |M\Delta_{n+1} - s|$. So this implies that $M\Delta_{n+1-w(\beta)} \prec M\Delta_{n+1} - s$ by Theorem 3.1.vi.

Then if we let $g$ be the generator of $M\Delta_{n+1-w(\beta)}$ and $h$ be the generator of $M\Delta_{n+1} - s$, then by Theorem 3.1.v, $g < h$ ($g$ cannot equal $h$ since the two delta-sets have a different number of elements). Thus we can add $s$ to each generator and obtain $\overline{M\Delta_{n+1-w(\beta)} + s} \prec M\Delta_{n+1}$ by Theorem 3.1.v again. Therefore, $\beta - s \prec M\Delta_{n+1}$.

At this point we simply apply Lemma 5.1 to reverse the initial push upwards and achieve $\beta \prec \overline{M\Delta_{n+1} + s}$. If we let $A = \beta \cap \{1, 2, ..., s\}$ and $B = \overline{M\Delta_{n+1} + s} \cap \{1, 2, ..., s\}$, then Lemma 5.1 applies because $w(\beta) > w(M\Delta_{n+1})$ implies that $|A| = w(\beta) \geq w(M\Delta_{n+1}) + 1 = w(\overline{M\Delta_{n+1} + s}) = |B|$.

Second case: $h(\beta) > h(M\Delta_{n+1})$

The proof of this case is almost identical to the previous one, except here we consider $\beta - t$ and $M\Delta_{n+1} - t$.

Third case: $w(\beta) = w(M\Delta_{n+1}) < s$ and $h(\beta) = h(M\Delta_{n+1})$

Again, consider the delta-set $D$ defined as the minimal delta-set with respect to partition containment which has $\beta$ as a subset. By Lemma 5.5, $D$ has the same height and width as $\beta$ in this case. Since the width of $D$ and $M\Delta_{n+1}$ are both smaller than $s$, the elements of them can be expressed uniquely as a difference of their respective generators and multiples of $s,t$.

We then want to compare the heights of the columns of $D$ and of $M\Delta_{n+1}$. We call the $i$th column of a delta-set to be the $i$th column to the right of the generator. The 0th column of any delta-set with unique element representation must contain a number of elements equal to its height, because for each element $x$ in the 0th column we can subtract $t$ until we achieve its residue modulo $t$, $x - b_x t$, where $0 < x - b_x t \leq t$. These values are distinct, because we have unique element representation, and furthermore no other residues are possible because they are not elements of the delta-set. Thus the number of elements in the 0th column equals the height of that delta-set.

Therefore, $D$ and $M\Delta_{n+1}$ must contain the same number of elements since both have the same height as $\beta$. From this we can deduce that the $i$th column of $D$ can be at most one taller than the $i$th column of $M\Delta_{n+1}$, as we will now show.

Suppose that instead the $i$th column of $D$ was two or more taller than the $i$th column of $M\Delta_{n+1}$. Suppose, more precisely, that the $i$th column of $D$ has $l$ elements and the $i$th column of $M\Delta_{n+1}$ has at least $l+2$, then the bottom element of the column in $D$ equals $y_1 + (l-1)s$, where $y_1$ is some value between 1 and $s$, and the bottom element of the column in $M\Delta_{n+1}$ is at least $y_2 + (l+1)s$, where $y_2$ is some value between 1 and $s$. Thus the bottom element of the $i$th column of $D$ must be at least $s$ more than the bottom element of the $i$th column of $M\Delta_{n+1}$.

But the bottom element of the $(i-1)$ st column of $D$ then is $y_1 - t + (l-1)s$, and the bottom element of the $(i-1)$ st column of $M\Delta_{n+1}$ is $y_2 - t + (l-1)s$ and their difference is again at least $s$. Thus even in the 0th column the largest elements are still at least $s$ apart, but that in turn means that by the same argument, the 0th column of $D$ has more elements than the 0th column of $M\Delta_{n+1}$, which is a contradiction.

Thus each column of $D$ can at most be one taller than the corresponding column of $M\Delta_{n+1}$ (except for the leftmost columns, which must be equal in size), and so we have that the difference in cardinalities between $D$ and $M\Delta_{n+1}$ must be at most one for every column besides the 0th. Hence, $|D| - |\beta| = |D| - |M\Delta_{n+1}| < w(D)$.

Therefore $\beta \prec D$ by Lemma 5.2. Since $\left|\overline{M\Delta_{n+1} + s}\right| - |M\Delta_{n+1}| \geq w(M\Delta_{n+1})$ we have that $D \prec \overline{M\Delta_{n+1} + s}$.

Final case: $w(\beta) = w(M\Delta_{n+1}) = s$

In this case we note that 1 is an element of both $\beta$ and $M\Delta_{n+1}$ since both beta-sets contain elements from every column and are both $(s,t)$-closed. Therefore we have that $P(\beta - 1)_i = P(\beta)_i$ and $P(M\Delta_{n+1})_i = P(M\Delta_{n+1} - 1)_i$, for $1 \leq i \leq n$, and $P(\beta)_{n+1} =$ $P(M\Delta_{n+1})_{n+1} = 1$. Thus anything we know about the relative values of the partitions for $\beta - 1$ and $M\Delta_{n+1} - 1$ will also be true about the relative values of $\beta$ and $M\Delta_{n+1}$.

So instead of comparing $\beta$ and $M\Delta_{n+1}$ we can compare $\beta - 1$ and $M\Delta_{n+1} - 1$ and see which of the cases above these new beta-sets fit into. If they fit into this final case again, we just subtract 1 again and again. Since $\beta$ has a finite number of elements, it would be impossible to continue doing this forever and ever, at some point, say $k$, we must have that $w(\beta - k) < w(\beta)$. If we let $k$ be the smallest value for which that is true then we must have that $w(\beta - k) = s - 1$, since $w(\beta - k + 1) = s$ by construction of $k$ and since subtracting 1 can remove at most one element (and hence one residue class modulo $s$).

But we also note that since subtracting 1 from $M\Delta_{n+1}$ causes it to shrink in size $M\Delta_{n+1} - 1 = M\Delta_n$, therefore we are eventually comparing $\beta - k$ and $M\Delta_{n+1-k}$ under one of the other cases. But we have proved in all the other cases that $\beta - k \prec \overline{M\Delta_{n+1-k} + s}$ and this latter set by the same reasoning is just $M\Delta_{n+1-k+s}$.

Then we can add $k$ back to our two sets to find that $\beta \prec M\Delta_{n+1+s} = \overline{M\Delta_{n+1} + s}$ as desired.

Q.E.D.

## 6. Further Investigations

Since Anderson, Olsson, and Stanton among others extended the study of partitions which are $t$-core to study partitions which are simultaneously $s$- and $t$-core, it makes sense to further extend the study and ask what we know about partitions which are simultaneously $t_1, t_2, \ldots, t_n$-core. We say that a partition is $T$-core where $T = \{t_1, t_2, \ldots, t_n\}$, if the partition is $t_i$-core for $1 \le i \le n$. Again, by convention we want $t_1 < t_2 < \cdots < t_n$ and the $t_i$ to be pairwise relatively prime.

Most of the results proved in this paper can be extended to the $T$-core case using the same method of proof, with two notable exceptions.

First, the maximal theorem itself is no longer true in the $T$-core case. Note that in the $(s,t)$-core case, being $(s,t)$-core is equivalent to being $(s,t)$-closed and set-wise contained in $\beta_{s,t} = \overline{\{st - s - t\}}$. However, in the $T$-core case, being $T$-core is equivalent to being $T$-closed and set-wise contained in $\beta_T$, but $\beta_T$ is no longer necessarily a delta-set. For example, in the $(5,6,7)$-core case, $\beta_T = \{1, 2, 3, 4, 8, 9\}$. (For information on how to derive $\beta_T$, see Anderson [1].)

Second, the canonical form must be defined differently. If $T$ has only two elements, then we use our old definition of $D = \{m - a_1 t_1 - a_2 t_2 \mid a_1 + a_2 \le b\}$ for appropriate choices of $m$ and $b$. Here, $D$ is 2-dimensional. We also have diagonals of the form $D_i = \{m - a_1 t_1 - a_2 t_2 \mid a_1 + a_2 = i\}$, which are in essence 1-dimensional, since if we choose our value for $a_1$ then we have only one choice for the value of $a_2$ and hence only one possibility for $m - a_1 t_1 - a_2 t_2$ .

However, when we extend $D$ to having more than 2 dimensions, the diagonals are quite different. If $D = \{ m - a_1 t_1 - \cdots - a_n t_n \mid a_1 + \cdots + a_n \leq b \}$, then the diagonals are $D_i = \{ m - a_1 t_1 - \cdots - a_n t_n \mid a_1 + \cdots + a_n = i \}$, which are $(n-1)$-dimensional. At first there doesn't appear to be a simple way of using these diagonals to make canonical forms until we realize that the diagonals can be rewritten, as $D_i = \{ m - a_1 t_1 - \cdots - a_n t_n \mid a_1 + \cdots + a_n = i \}$ $= \{ (m - i t_1) - a_2'(t_2 - t_1) - \cdots - a_n'(t_n - t_1) \mid a_2' + \cdots + a_n' \leq i \}$ . Which is just another lattice of the same form as $D$ in one less dimension. So if we know how to construct a canonical form in $(n-1)$ dimensions, then we can begin to construct a notion of canonical forms in $n$ dimensions. Then we can generalize the first step in constructing the canonical form by saying that if $D_i$ has $n$ elements in it, then we replace it by the canonical beta-set with n elements in $D_i$. Similarly we generalize the second step by saying that if we can reduce the canonical form in $D_i$ by 1 element and increase the canonical form in $D_{i-a}$ by 1 element while maintaining $T$-closure in $D$, then we do so.

It is still an open question as to whether the error term of the second generalization (both in the $(s,t)$-core and $T$-core case) can be shrunk. We conjecture, at the moment, that the term can be shrunk to $t_1 - w_1(\beta)$, where $w_1(\beta)$ equals the number of distinct residue classes of $\beta$ modulo $t_1$.

**Glossary of Mathematical terms and symbols**

$\in$ – Denotes "is an element of". $a \in A$ means that $a$ is an element of the set $A$

$\subset$ – Denotes "is a subset of". $A \subset B$ means that every element in $A$ is also an element of $B$.

$\cup$ – Denotes "union". $A \cup B$ is the set of all elements that are in $A$ or $B$ (or both)

$\cap$ – Denotes "intersection". $A \cap B$ is the set of all elements that are in *both* $A$ and $B$.

$|\ |$ – Denotes "cardinality". $|A|$ is the number of elements in the set $A$.

$:=$ – Denotes "definition". $E_x := \sin x$ means that we define $E_x$ by $\sin x$.

$\prec$ – Denotes "partition-wise containment" for beta-sets. See section 2.

$<$ – Denotes "partition-wise containment" for partitions. See section 2.

$\{\ |\ \}$ – Denotes a set with restrictions. To the left of the bar is some formula featuring variables, and on the right of the bar are the restrictions of those variables. The set then is all possible elements that can be generated under those restrictions.

$\bigcup$ – Also denotes union. This is used as shorthand for writing down several union

symbols. $\displaystyle\bigcup_{k=1}^{n} A_k = A_1 \cup A_2 \cup \cdots \cup A_n$

\ – Denotes subtraction of sets. $A\backslash B$ is the set of all elements in $A$ that are not in $B$.

$\Sigma$ – Denotes a shorthand for summation. $\displaystyle\sum_{k=1}^{n} a_k = a_1 + a_2 + \cdots + a_n$

$\mathbb{N}$ – Denotes the set of natural numbers. $\mathbb{N} = \{1, 2, 3, \ldots\}$

$\mathbb{Z}$ – Denotes the set of integers. $\mathbb{Z} = \{\ldots, -2, -1, 0, 1, 2, \ldots\}$

Bijection – A function $f$ is a bijection if $f$ is one-to-one ( $f(x) = f(y)$ implies $x = y$ ) and onto (for every $y$ there exists some $x$ such that $f(x) = y$ ).

Lattice – for the purposes of this paper, a lattice is just a subset of the integer coordinate points in $n$-dimensional space.

Matrix multiplication – Given two matrices $A = \begin{pmatrix} a_{11} & \cdots & a_{1n} \\ \vdots & \ddots & \vdots \\ a_{m1} & \cdots & a_{mn} \end{pmatrix}$, $B = \begin{pmatrix} b_{11} & \cdots & b_{1n} \\ \vdots & \ddots & \vdots \\ b_{m1} & \cdots & b_{mn} \end{pmatrix}$, we have that $AB = \begin{pmatrix} c_{11} & \cdots & c_{1n} \\ \vdots & \ddots & \vdots \\ c_{m1} & \cdots & c_{mn} \end{pmatrix}$ where $c_{i,j} = \sum_{k=1}^{n} a_{i,k} b_{k,j}$ .

Modulo – Two numbers are said to be equivalent modulo $n$ if their remainders when divided by $n$ are equal. For example, 5 and 9 are equivalent modulo 4 (both have a remainder of 1).

Triangular number – the $n$th triangular number is the some of the first $n$ integers. Thus the first few triangular numbers are $1 = 1,\ 1 + 2 = 3,\ 1 + 2 + 3 = 6,\ 1 + 2 + 3 + 4 = 10$ .

Tuple – A tuple is a list or ordered set of elements. An $n$-tuple is a tuple containing $n$ elements.

Relatively prime – two numbers are relatively prime if there are no prime numbers that divide both evenly.